\documentclass[a4paper,10pt]{article}
\usepackage{graphicx,multirow}
\usepackage{amssymb,amsmath,amsthm}
\usepackage{lscape,lipsum,microtype}
\usepackage[margin=2.54cm]{geometry}  
\usepackage{hyperref}

\usepackage{ifluatex}
\ifluatex
 \usepackage{unicode-math}
\fi

\usepackage{latexsym}
\usepackage{amsmath}
\usepackage{amsfonts}
\usepackage{amssymb}
\usepackage{amscd}

\renewcommand{\d}{\delta }
\newcommand{\D }{\Delta }

\renewcommand{\l }{\lambda }

\newcommand{\n }{\nabla }

\newcommand{\intbar}{\mathop{\int\makebox(-13.5,0){\rule[4pt]{.7em}{0.3pt}}%
\kern-6pt}\nolimits}

\newcommand{\be}{\begin{equation}}
\newcommand{\ee}{\end{equation}}
\newcommand{\bes}{\begin{equation*}}
\newcommand{\ees}{\end{equation*}}
\newcommand{\ba}{\begin{eqnarray}}
\newcommand{\ea}{\end{eqnarray}}
\newcommand{\bas}{\begin{eqnarray*}}
\newcommand{\eas}{\end{eqnarray*}}
\newenvironment{pf}{\noindent{\sc Proof}.\enspace}{\rule{2mm}{2mm}\medskip}
\newenvironment{pfn}{\noindent{\sc Proof}}{\rule{2mm}{2mm}\medskip}

\newcommand{\R}{\mathbb{R}}

\newcommand{\Z}{\mathbb{Z}}

\newcommand{\N}{\mathbb{N}}

\author{ Mohammed ALDAWOOD$^a$, Cheikh Birahim NDIAYE$^b$}

\date{}

\title{\bf Cherrier-Escobar problem for the elliptic Schr\"odinger-to-Neumann map.}
\begin{document}

\newtheorem{lem}{Lemma}[section]
\newtheorem{pro}[lem]{Proposition}
\newtheorem{thm}[lem]{Theorem}
\newtheorem{rem}[lem]{Remark}
\newtheorem{cor}[lem]{Corollary}
\newtheorem{df}[lem]{Definition}

\maketitle

\begin{center}
{\small

\noindent  $^{a, b}$\; Department of Mathematics Howard University \\  Annex 3, Graduate School of Arts and Sciences\\ DC 20059 Washington, USA.
}

\end{center}

\footnotetext[1]{E-mail addresses: cheikh.ndiaye@howard.edu, mohammed.aldawood@bison.howard.edu\\
\thanks{\\ C. B. Ndiaye was partially supported by NSF grant DMS--2000164.}}

\
\

\begin{center}
{\bf Abstract}
\end{center}
In this paper, we study a Cherrier-Escobar problem for the extended problem corresponding to the elliptic Schr\"odinger-to-Neumann map on a compact $3$-dimensional Riemannian manifold with boundary. Using the algebraic topological argument of Bahri-Coron\cite{bc}, we show solvability under the assumption that the extended problem corresponding to the elliptic Schr\"odinger-to-Neumann map has a positive first eigenvalue, a positive Green's function, and also verifies the strong maximum principle.
 \begin{center}


\bigskip

\centerline{\bf AMS subject classification: 53C21, 35C60, 58J60, 55N10.}

\end{center}

\section{Introduction and statement of the results}
Boundary Value problems (BVP) with critical nonlinearity on the boundary of the form    
\begin{equation}\label{BVP1}
\left\{
\begin{split}
-\D_g u+qu&=0&\;\;\;\;\text{in}\;\;\;\;&M,\\
-\frac{\partial u}{\partial n_{g}}+\frac{n-2}{2(n-1)}H_{g}&=u^{\frac{n}{n-2}}&\;\;\;\;\text{on}\;\;\;\;&\partial M,\\
u&>0\;\;\;\;&\text{on}\;\;\;\;&\overline{M},
\end{split}
\right.
\end{equation}
have received a lot of attention in the last decades. In \eqref{BVP1}, \;$\left(\overline{M},g\right)$\; is a $n$-dimensional compact Riemannian manifold with boundary \;$\partial M$\; and interior \;$M,$\; with \;$n\geq3.$\; Furthermore, \;$\Delta_g$\; denotes the Laplace–Beltrami with respect to \;$g,$\;\;$\frac{\partial}{\partial n_{g}}$\; denotes the inner Neumann derivative with respect to\;$g,$\;\;$H_g$\; denotes the mean curvature of \;$\partial M$\; with respect to\;$g$\; in the normal direction, and the potential \;$q$\; is a bounded smooth function defined on\;$M.$\; In this paper,
we study the particular case \;$dim\;\overline{M}=3$\; and \;$H_g=0.$\; Hence, the BVP of interest becomes 
\begin{equation}\label{BVP2}
\left\{
\begin{split}
-\D_g u+qu&=0&\;\;\;\;\text{in}\;\;\;\;&M,\\
-\frac{\partial u}{\partial n_{g}}&=u^{3}&\;\;\;\;\text{on}\;\;\;\;&\partial M,\\
u&>0&\;\;\;\;\text{on}\;\;\;\;&\overline{M}.
\end{split}
\right.
\end{equation}
\vspace{4pt}

\noindent
Looking at the Cherrier-Escobar problem studied in \cite{martndia2} as the $\frac{1}{2}$-Yamabe problem, we have BVP \eqref{BVP2} is related to a Cherrier-Escobar problem for the elliptic Schr\"odinger-to-Neumann map \;$P_q$\;  defined by \(\ P_q:{C^{\infty}(\partial M)} \to {C^{\infty}(\partial M)}\), \;\(u\to \frac{-\partial U_q}{\partial n_{g}}\)  where \;$U_q$\; is the unique solution of
\begin{equation}\label{BVP3}
\left\{
\begin{split}
-\D_g U_q+qU_q&=0&\;\;\;\;\text{in}\;\;\;\;&M,\\
U_q&=u&\;\;\;\;\text{on}\;\;\;\;&\partial M.
\end{split}
\right.
\end{equation}
Indeed, \;$u$\; is solution of BVP \eqref{BVP2} implies 
\begin{equation}\label{PQ}
\left\{
\begin{split}
 P_qu&=u^{3}&\;\;\;\;\text{on}\;\;\;\;&\partial M,\\
u&>0&\;\;\;\;\text{on}\;\;\;\;&\partial M,
\end{split}
\right.
\end{equation}
and \;$u$\; is solution of \eqref{PQ} implies that \;$U_q$\; defined by \eqref{BVP3} is solution of BVP \eqref{BVP2} thanks to fact that the extended problem associated to \;$P_q$\; verifies the strong maximum principle as we will assume (see \eqref{pq2}-\eqref{pq3} below). We used the terminology the elliptic Schr\"odinger-to-Neumann map, since \;$P_q$\; is clearly the Dirichlet to Neumann associated to Schr\"odinger operator \;$-\Delta_{g}+q.$\;
\vspace{4pt}

\noindent
As in \cite{aldawood}, it is easy to see that a necessary condition of the existence of solution to \eqref{BVP2} is that the first eigenvalue \;$\lambda_{1}(P_q)$\; of the extended problem corresponding to \;$P_q$\; i.e. of the eigenvalue problem 
\begin{equation}\label{BVP4}
    \left\{
\begin{split}
-\D_g u+qu&=0&\;\;\;\;\text{in}\;\;\;\;&M,\\
-\frac{\partial u}{\partial n_{g}}&=\lambda u&\;\;\;\;\text{on}\;\;\;\;&\partial M,
\end{split}
\right.
\end{equation}
is strictly positive.\\
\vspace{4pt}

\noindent
We use the symbol \;$\lambda_{1}(P_q),$\; since it can be easily checked that it corresponds also to the first eigenvalue of \;$P_q.$\; Hence, from now on we will assume that
\begin{equation}\label{pq1}
    \lambda_1(P_q)>0.
\end{equation}
\vspace{2pt}

\noindent
As already said, we will assume that the extended problem corresponding to \;$(P_q)$\; satisfies the strong maximum principle, in the sense that 
\begin{equation}\label{pq2}
\left\{
\begin{split}
P_{q}u&\geq0\;\;\;\;&\text{on}\;\;\;\;\partial M,\\
u&\geq0\;\;\;\;&\text{on}\;\;\;\;\partial M,\\
u&\neq0\;\;\;\;&\text{on}\;\;\;\;\partial M,
\end{split}
\right.
\end{equation}
implies 
\begin{equation}\label{pq3}
    U_{q}>0\;\;\;\;\;\;\text{on}\;\;\;\;\;\overline{M}
\end{equation}
with \;$U_{q}$\; defined by \eqref{BVP3}. Thus thanks to elliptic regularity theory (see \cite{pc}), a smooth solution of \eqref{BVP2} can be found by looking at critical points of the functional \;$J_q$\; defined by 
\begin{equation}\label{eq:dfjq}
J_q(u):=\frac{\left<u, u\right>_q}{(\oint_{\partial M}u^4\;dS_g)^{\frac{1}{2}}}, \;\;\;\;u\in H^{1}_+(M):=\{u\in H^1(M):\;u\geq 0 \;\;\text{and}\;\;u\neq 0\},
\end{equation}
with 
\begin{equation}\label{scalq}
\left<u, u\right>_q=\int_{M}\left(|\n_g u|^{2}+qu^{2}\right)\;dV_g,
\end{equation}
\vspace{2pt}

\noindent
where \;$\nabla_g$\; denotes the covariant derivative with respect to \;$g,$\;\;$dV_{g}$\; denotes the volume form with respect to \;$g,$\; and \;$dS_{g}$\; denotes the volume form with repect to the Riemannian metric induced by \;$g$\; on \;$\partial M.$\; Moreover, in formula \eqref{eq:dfjq}, \;$H^1(M)$\; denotes the usual Sobolev space of functions which are \;$L^2$-integrable together with their first derivatives.
\vspace{4pt}

\noindent
The existence of solutions for \eqref{BVP2} are known under the assumption \;$q=R_{g}$\; up to a positive constant, where \;$R_{g}$\; denotes the scalar curvature of \;$\left(\overline{M},g\right).$\; There are a lot of works related to equation \eqref{BVP1}, see \cite{sma1},\cite{sma2}, \cite{fcm1}, \cite{fcm2}, \cite{jfe1}, \cite{jfe2}, and \cite{moa}. In this paper, we use the Barycenter technique of Bahri-Coron\cite{bc} to study \eqref{BVP2} for general potentials \;$q.$\; We prove the following result:
\noindent
\begin{thm}\label{thm1}
Let \;$\left(\overline{M},g\right)$\; be \;$3$-dimensional compact Riemannian manifold with boundary \;$\partial M$\; and interior \;$M$\; such that \;$H_g=0.$\; Let also \;$q$\; be a bounded smooth potential defined on \;$M.$\; Assuming that \;$\lambda_{1}(P_q)>0$\; (see \eqref{BVP4} for the definition), the extended problem corresponding to \;$P_q$\; satisfies the strong maximum principle (in the sense of \eqref{pq2}-\eqref{pq3}), and the Green's function \;$G$\; of the extended problem corresponding to \;$P_q$\; defined by \eqref{eqgreen} is strictly positive, then the BVP \eqref{BVP2} has a least one solution.
\end{thm}
\vspace{4pt}

\noindent
As in \cite{aldawood}, to prove Theorem \ref{thm1}, we will use the Barycenter technique of Bahri-Coron\cite{bc} which is possible since \;$dim\;\overline{M}=3$\; and \;$H_g=0$\; imply the problem under study is a Global one (for the definition of "Gobal" for Yamabe and Cherrier-Escobar type problems, see \cite{nss}). Indeed, as in \cite{aldawood} and \cite{nss}, we will follow the scheme of the Algebraic topological argument of Bahri-Coron\cite{bc} as performed in the work \cite{martndia2} of the second author and Mayer. As in \cite{aldawood}, one of the main difficulty with respect to the works  \cite{martndia2} and \cite{nss} is the presence of the linear term \;$"qu"$\; and the lack of conformal invariance. As in \cite{aldawood}, to deal with such a difficulty, we use the fact that \;$dim\;\overline{M}=3,$\;\;$H_g=0,$\; and the Brendle\cite{bre1}-Schoen\cite{sc}' s bubble construction to run the scheme of the Algebraic topological argument of Bahri-Coron\cite{bc} for  the existence.
\vspace{4pt}

\noindent
The plan of the paper is as follows. In Section \ref{NP}, we fix some notations and discuss some preliminaries. In Section \ref{PSS}, we recall the profile decomposition of Palais-Smale (PS)-sequences for \;$J_q.$\; We also introduce the neighborhoods of potential critical points at infinity of \;$J_q$\; and  their associated selection maps. We also state a Deformation Lemma for \;$J_q$\; taking into account the possible bubbling phenomena involved in the study \eqref{BVP2}. In Section \ref{SAE}, we derive some sharp self-action estimates needed for the application of the Barycenter technique of Bahri-Coron\cite{bc} for existence. In Section \ref{IE}, we derive some sharp inter-action estimates needed for the application of the Barycenter technique of Bahri-Coron\cite{bc} for existence. In Section \ref{ATA},  we present the algebraic topological argument for existence. In Section \ref{APP}, we collect some technical estimates. 
%
%
%
%
%
%
\section{Notations and preliminaries}\label{NP}
In this Section, we fix some notations and discuss some preliminaries. We start by introducing some notations.\\
\vspace{4pt}

\noindent
For \;$x=\left(x_1,x_2,x_3\right)\in \bar{\R}^3_{+}=\R^2\times\bar{\R}_{+},$\; with \;$\bar{\R}_{+}=[0,+\infty),$\; we set \;$\bar{x}=(x_1,x_2)\in \R^{2},$\; so that \;$x=\left(\bar{x},x_3\right),$\; and we set also \;$r=|x|=\sqrt{x^{2}_1+x^{2}_2+x^{2}_3}.$\; We define \;$\R^3_{+}=\R^2\times\R_{+}$\; with \;$\R_{+}=(0,\infty).$\; We identify the boundary of \;$\R^3_{+}$\; denoted \;$\partial \R^3_{+}$\; with \;$\R^{2}.$\; For \;$r>0,$\; we define the half ball \;$B^{+}_{r}(0)$\; with radius\;$r$\;centered at \;$0,$\; by $$B^{+}_{r}(0)=B_{r}(0)\cap \bar{\R}^3_{+},$$ 
with \;$B_{r}(0)$\; denoting the Euclidean ball centered at \;$0$\; with radius \;$r$\; in \;$\R^{3},$\; namely
$$B_{r}(0)=\{x\in \R^3:\; |x|<r \}.$$ 
Furthermore, we define 
$$\partial B^{+}_{r}(0)=B_{r}(0)\cap\R^{2}.$$\\
Large positive constants are usually denoted by \;$C$\; and the value of \;$C$\; is allowed to verify from formula to formula and also within the same line.\\\\  
For \;$r>0$\; and \;$a\in \partial M,$\; we define
\begin{equation*}
    \begin{split}
    &B(a,r)=\{x\in M:\;d_g(a,x)<r\},\\
     &\hat{B}(a,r)=\{x\in \partial M:\; d_g(a,x)<r\},\\
     &\hat{B}_{r}(0)=\{x\in \R^{2}:\; |x|<r\}.
    \end{split}
\end{equation*}
We set 
\begin{equation}\label{uqsqrt}
    \|u\|_{q}=\sqrt{\left<u, u\right>_q}\;,\;\;u\in H^{1}(M),
\end{equation}
and \;$\left<u, u\right>_q$\; is as in \eqref{scalq}.\\ 
\vspace{4pt}

\noindent
For \;$\lambda>0,$\; we set 
\begin{equation}\label{bubble-M3}
     \delta_{0,\lambda}(x)=c_0\left[\frac{\lambda}{\left(1+\lambda x_3\right)^2+\lambda^2\left|\bar{x}\right|^2}\right]^{\frac{1}{2}},\;\;\;\;x=(\bar{x},x_3)\in \bar{\R}^3_{+}
\end{equation}
with \;$c_0>0$\; such that
\vspace{4pt}

\noindent
\begin{equation*}
\left\{
\begin{split}
\D \delta_{0,\lambda}&=0&\;\;\;\;\text{in}\;\;\;\;&\R^3_{+},\\
-\frac{\partial \delta_{0,\lambda}}{\partial x_{3}}&=\delta^{3}_{0,\lambda}\;\;\;\;&\text{on}\;\;\;\;&\partial\R^3_{+},
\end{split}
\right.
\end{equation*}
where
$$\Delta=\sum_{i=1}^3 \frac{\partial^{2}}{\partial x^{2}_{i}}$$ is the Euclidean Laplacian on \;$\R^{3}.$\;\\
\vspace{2pt}

\noindent
It is a well-know fact that (and easy to check) 
\begin{equation}\label{S3}
\int_{\R^3_{+}}|\nabla \delta_{0, \lambda}|^2=\int_{\R^3_{+}}|\nabla \delta_{0,1}|^2=\int_{\R^2}\delta_{0,\lambda}^4=\int_{\R^2}\delta_{0,1}^4
\end{equation}
Moreover, we set
\begin{equation}\label{S}
\mathcal{S}=\frac{\int_{\R^{3}_{+}}|\n \d_{0,1}|^2}{(\int_{\R^2}\delta_{0,1}^4)^{\frac{1}{2}}}, 
\end{equation}
We define 
\begin{equation}\label{c3}
c_1=\int_{\R^2}\left(\frac{1}{1+|y|^2}\right)^{\frac{3}{2}}dy.
\end{equation}
\vspace{4pt}

\noindent
For \;$a\in \partial M,$\; we let \;$G(a,x)$\; be the unique solution of
\begin{equation}\label{eqgreen}
\left\{
\begin{split}
-\Delta_g G(a,x)+qG(a,x)&=0,&\;\;\;&x \in M,\\
-\frac{\partial G(a,x)}{\partial n_{g}}&=2\pi\delta_a(x),&\;\;\;&x \in \partial M.
\end{split}
\right.
\end{equation}
\vspace{4pt}

\noindent
Since \;$q$\; is a bounded smooth function defined on \;$M,$\;\;$H_g=0,$\; and \;$dim\;\overline{M}=3,$\; then Green's function \;$G(a,x)$\; satisfies the following estimates 
  \begin{equation}\label{estg}
 \left|G(a,x)-\frac{1}{d_g(a,x)}\right|\leq C,\;\;\;\;\text{for}\;\;\;\;x\neq a\in\overline{M},
 \end{equation}
 and
 \begin{equation}\label{estgg}
 \left|\nabla \left(G(a,x)-\frac{1}{d_g(a,x)}\right)\right|\le \frac{C}{d_g(a,x)},\;\;\;\; \text{for}\;\;\;\;x\neq a\in \overline{M},
 \end{equation}
 with \;$C$\; a positive constant.\\
\vspace{2pt}

\noindent
Moreover, we have \;$\overline{M}$\; compact implies that there exists 
 \begin{equation}\label{delta0}
     \delta_0>0
 \end{equation}
  such that \;$\forall$\;\;$a\in \partial M$\; and \;$\forall$\;\;$0<2\delta\le \delta_0,$\; the Fermi coordinates centered at \;$a$\; defines a smooth map
 \begin{equation}\label{psi}
     \psi_a:B^+_{2\delta}(0)\to \overline{M},
 \end{equation}
 identifying a neighborhood \;$O(a)$\; of \;$a$\; in \;$\overline{M}.$\; We will identify a point \;$y=\psi_a(x)\in O(a)$\; with \;$x\in B^+_{2\delta}(0).$\; With this agreement and recalling that \;$H_g=0$\;, we have that an expansion of the Riemannian metric\;$g$\; and \;$\sqrt{|g|}$\; (where \;$|g|$\; denotes the modulus of the determinant of \;$g$\;) on \;$B^+_{2\delta}(0)$\; is given by the following formulas 
 \begin{equation}\label{E_Metric}
 \begin{split}
     g^{ij}(x)&=\delta_{ij}+2(L_g)_{ij}x_n(x)+\frac{1}{3}R_{ikjl}[\hat{g}]x_kx_l(x)+g^{ij}_{nk}x_nx_k(x)+\left\{3(L_g)_{ij}(L_g)_{kj}+R_{injn}[g]\right\}x^{2}_n(x)\\
     &+o(|x|^3),\;\;\;\;\;x\in B^+_{2\delta}(0)\\
     \sqrt{|g(x)|}&=1-\frac{1}{6}R_{ic}[\hat{g}]_{ij}x_i x_j(x)-\left[\frac{1}{2}\Arrowvert L_g\Arrowvert^2+R_{ij}[g]_{nn}\right]x^2_{n}(x)+o(|x|^3),\;\;\;\;\;x\in B^+_{2\delta}(0).
     \end{split}
 \end{equation}
 In the formulas in \eqref{E_Metric}, \;$n=3,$\;\;$(L_g)_{ij}$\; denotes the component of the second fundamental form of \;$\partial M$\; with respect to \;$g,$\; \;$(R_{ic}[g])_{ab},$ $a,b=1,..,n$\; denotes the component of the Ricci tensor of \;$\overline{M}$\; with respect to \;$g,$\;\;$\hat{g}:=g|_{\partial M}$\; is the Riemannian metric induced by \;$g$\; on \;$\partial M,$\;\;$(R_{ic}[\hat{g}])_{ij},$\;\;$i,j=1,..,n-1$\; are the components of the Ricci tensor of \;$\partial M$\; with respect to \;$\hat{g},$\;\;$R_{abcd}[g],$\;\;$a,b,c,d=1,..,n$\; denotes the components of Riemann tensor of \;$\overline{M}$\; with respect to \;$g,$\; and \;$R_{ijkl}[\hat{g}],$\;\;$i,k=1,..,n-1$\; denotes the components of Riemann curvature tensor of \;$\partial M$\; with respect to \;$\hat{g}.$\; All the tensors in the right of \eqref{E_Metric} are evaluated at \;$0,$\;and we also use Einstein summation convention for repeated indexes.
 \vspace{4pt}
 
 \noindent
 Let \(\ \chi:\mathbb{R} \to\R\)\; be a smooth cut-off function satisfying
\begin{equation}\label{chi}
\chi(t)= \left\{ \begin{array}{ll}
         1,&\;\;\mbox{if\;\;\;$t\leq  1$},\\
        0,&\;\;\mbox{if\;\;\;$t\geq 2$}.\end{array} \right.
        \end{equation}
\vspace{6pt}
 
 \noindent
For \;$0<2\delta<\delta_0,$\; we define
\begin{equation}\label{chid}  
        \chi_{\delta}(x)=\chi\left(\frac{|x|}{\delta}\right),\;\;\;\;x\in \bar{\R}^3_{+}.
        \end{equation}
\vspace{6pt}
 
 \noindent
For \;$0<2\delta<\delta_0,$\; \;$a\in \partial M,$\; and \;$\lambda>0,$\; we define the Brendle\cite{bre1}-Schoen\cite{sc}'s bubble 
\begin{equation}\label{uald}
       u_{a, \lambda}(x)=u_{a,\lambda,\delta}(x):=\chi^{a}_\delta(x)\hat{\delta}_{a,\lambda}(x) +(1-\chi^{a}_\delta(x))\frac{c_0}{\sqrt{\lambda}}G(a,x),\;\;\;\;\text{for}\;\;\;\;x\in \overline{M},
 \end{equation}
 \noindent
 where 
 \begin{equation}\label{deltahat}
  \hat{\delta}_{a,\lambda}(x)=\delta_{0,\lambda}\left(\psi^{-1}_{a}(x)\right),
 \end{equation}
 and
 \begin{equation}\label{chi1}
  \chi^{a}_\delta=\chi_\delta(\psi^{-1}_a(x)),
 \end{equation}
 with
 \begin{equation}\label{psi0}
  \psi_a:B^+_{\delta_0}(0)\to \overline{M}.
 \end{equation}
 \vspace{4pt}
 
 \noindent
 Thus, recalling that we are under the assumption \;$G>0$\; (for the definition of \;$G$\; see \eqref{eqgreen}), then \;$\forall$\;\;$a\in \partial M$\; and \;$\forall$\;\;$0<2\delta<\delta_0,$\; we have 
\begin{equation}\label{ual}
u_{a,\lambda}\in H^1(M), \;\;\;\;\text{and}\;\;\;\;u_{a, \l}>0\;\;\;\;\text{in}\;\;\;\;\overline{M}.
\end{equation}
For \;$a_i, a_j\in \partial M,$\; and \;$\l_i, \l_j>0,$ we define
\begin{equation}\label{varepij}
\varepsilon_{ij}=\left[\frac{1}{\frac{\lambda_i}{\lambda_j}+\frac{\lambda_j}{\lambda_i}+\lambda_i\lambda_jG^{-2}(a_i,a_j)}\right]^{\frac{1}{2}}.
\end{equation}
\vspace{4pt}
 
 \noindent
Moreover,\;\;for \;$0<2\delta<\delta_0,$\;\;$a_i, a_j\in \partial M,$\; and \;$\l_i, \l_j>0,$\; we define
\begin{equation}\label{epij}
\epsilon_{ij}=\oint_{\partial M} u^3_{a_i,\lambda_i}u_{a_j,\lambda_j}\;dS_g
\end{equation}
and
\begin{equation}\label{eij}
e_{ij}=\int_{M}\nabla_g u_{a_i,\lambda_i}\nabla_g u_{a_j,\lambda_j}dV_g+\int_{M}q u_{a_i,\lambda_i}u_{a_j,\lambda_j}\;dV_g.
\end{equation}
\vspace{4pt}
 
 \noindent
 To end the section, we derive the following $C^{0}$-estimate needed for the energy and the inter-action estimates required for the application of the Barycenter technique of Bahri-Coron\cite{bc} for existence.
 %
 %
 %
 %
 %
 %
 \begin{lem}\label{c0estimate}
 Assuming that \;$\theta>0$\; is small, then there exists \;$C>0$\; such that \;$\forall$\;\;$a \in \partial M,$\;\;$\forall$\;\;$0<2\delta<\delta_0$\; and \;$\forall$\;\;$0<\frac{1}{\l}\le\theta\delta,$\; we have
\begin{equation}\label{C_01}
\begin{split}
    \left|\left(-\Delta_g +q\right)u_{a,\lambda}(x)\right|\le C&\left[\frac{1}{\delta^2\sqrt{\lambda}}\large{1}_{\{y\in M:\;\delta\le d_g(a,y)\le 2\delta\}}(x)+\hat{\delta}_{a,\lambda}(x) \large{1}_{\{y\in M:\;d_g(a,y)\le 2\delta\}}(x)\right.\\
  +&\left.\quad\left(\frac{\lambda}{1+\lambda^{2}d^{2}_g(a,x)}\right)^{\frac{1}{2}}\large{1}_{\{y\in M:\;d_g(a,y)\le 2\delta\}}(x)\right],\;\;\;\forall \;\;x\in M,
  \end{split}
\end{equation}
and
\begin{equation}\label{C_02}
    \left|-\frac{\partial u_{a,\lambda}(x)}{\partial n_{g}} -u^3_{a,\lambda}(x)\right|\le  C\left[\left(\frac{\lambda}{1+\lambda^{2}d^{2}_g(a,x)}\right)^{\frac{3}{2}}\Large{1}_{\{y\in \partial M:\;d_g(a,y)\geq \delta\}}(x)\right],\;\;\;\forall \;\;x\in \partial M.
\end{equation}
\end{lem}
 
\begin{pf}
To simplifies notation, let us set \;  $G_a(\cdot):=G(a,\cdot)$\; and \;$\bar G_a(\cdot)=c_0G_a(\cdot)$. Then, we have 
\begin{equation}\label{CC_0}
u_{a,\lambda}(x)=\chi^{a}_\delta  (x)\hat{\delta}_{a,\lambda}(x)+(1-\chi^{a}_\delta(x))\frac{\bar{G}_a(x)}{\sqrt{\lambda}},\;\;\;\;x\in M.
\end{equation}
To deal with \eqref{C_01}, first we observe
\begin{equation*}
    \begin{split}
     \left(-\Delta_g+q\right) u_{a,\lambda}(x)&=\left(-\Delta_g+q\right)\left(\chi^{a}_\delta(x)\hat{\delta}_{a,\lambda}(x)+\left(1-\chi^{a}_\delta(x)\right)\frac{\bar{G}_a(x)}{\sqrt{\lambda}}\right)\\
     &=\left(-\Delta_g+q\right)\left[\chi^{a}_\delta(x)\left(\hat{\delta}_{a,\lambda}(x)-\frac{\bar{G}_a(x)}{\sqrt{\lambda}}\right)\right]+\left(-\Delta_g+q\right)\left[\frac{\bar{G}_a(x)}{\sqrt{\lambda}}\right],\;\;\;\;x\in M.   
    \end{split}
\end{equation*}
Now, since \;$x\in M$\; and \;$a\in \partial M,$\; then \;$\left(-\Delta_g+q\right)\left[\frac{\bar{G}_a(x)}{\sqrt{\lambda}}\right]=0.$\;
\vspace{4pt}
 
 \noindent
This implies
\begin{equation*}
    \begin{split}
     \left(-\Delta_g+q\right)u_{a,\lambda}(x)&=-\Delta_g\chi^{a}_\delta(x)\left[\hat{\delta}_{a,\lambda}(x)-\frac{\bar G_a(x)}{\sqrt{\lambda}}\right]-2\nabla_g\chi^{a}_\delta(x) \nabla_g\left[\hat{\delta}_{a,\lambda}(x)-\frac{\bar G_a(x)}{\sqrt{\lambda}}\right]\\
     &-\chi^{a}_\delta(x) \Delta_g\left[\hat{\delta}_{a,\lambda}(x)-\frac{\bar G_a(x)}{\sqrt{\lambda}}\right]+q\chi^{a}_\delta(x)\left[\hat{\delta}_{a,\lambda}(x)-\frac{\bar G_a(x)}{\sqrt{\lambda}}\right],\;\;\;\;x\in M.  
    \end{split}
\end{equation*}
Thus,
$$
\left(-\Delta_g+q\right)u_{a,\lambda}(x)=\sum_{m=1}^4 I_m
$$
with 
\begin{equation*}
\begin{split}
&I_1=-\Delta_g\chi^{a}_\delta(x)\left[\hat{\delta}_{a,\lambda}(x)-\frac{\bar G_a(x)}{\sqrt{\lambda}}\right],\;\;\;\;x\in M,\\
&I_2=-2\left<\nabla_g\chi^{a}_\delta(x), \nabla_g\left[\hat{\delta}_{a,\lambda}(x)-\frac{\bar G_a(x)}{\sqrt{\lambda}}\right]\right>,\;\;\;\;x\in M,\\
&I_3=q\chi^{a}_\delta(x)\hat{\delta}_{a,\lambda}(x),\;\;\;\;x\in M,\\
&I_4=-\chi^{a}_\delta(x) \Delta_g\hat{\delta}_{a,\lambda}(x),\;\;\;\;x\in M.
\end{split}
\end{equation*}
\vspace{4pt}

\noindent
We will estimate each of the \;$I_m$'s one at a time. For \;$I_1,$\; we have 
\begin{equation}\label{I_1}
I_1=-\Delta_g\chi^{a}_\delta(x)\left[\hat{\delta}_{a,\lambda}(x)-\frac{c_0}{\sqrt{\lambda}d_g(a,x)}+\frac{c_0}{\sqrt{\lambda}d_g(a,x)}-\frac{\bar G_a(x)}{\sqrt{\lambda}}\right],\;\;\;\;x\in M.
\end{equation} 
Applying \eqref{bubble-M3} and \eqref{estg}, we derive 
\begin{equation}\label{diff1}
\left|\hat{\delta}_{a,\lambda}(x)-\frac{c_0}{\sqrt{\lambda}d_g(a,x)}\right|\le \frac{C}{\sqrt{\lambda}},\;\;\;\;x\in M,
\end{equation}
and 
\begin{equation}\label{diff2}
\left|\frac{c_0}{\sqrt{\lambda}d_g(a,x)}-\frac{\bar G_a(x)}{\sqrt{\lambda}}\right|\le\frac{C}{\sqrt{\lambda}},\;\;\;\;x\in M. 
\end{equation}
Using \eqref{chi} and \eqref{chid}, we get

\begin{equation}\label{flchi}
|\Delta_g\chi^{a}_\delta(x)|\le \frac{C}{\delta^2}\large{1}_{\{y\in M:\;\delta\le d_g(a,y)\le 2\delta\}}(x),\;\;\;\;x\in M.
\end{equation}
Hence, combining  \eqref{I_1}-\eqref{flchi},  we obtain
\begin{equation}\label{estI1}
|I_1|\le \frac{C}{\delta^2\sqrt{\lambda}} \large{1}_{\{y\in M:\;\delta\le d_g(a,y)\le 2\delta\}}(x),\;\;\;\;x\in M.
\end{equation}
\vspace{6pt}
 
 \noindent
In order to estimate \;$I_2,$\; we first write
\begin{equation}\label{I_2}
I_2=2\left<\nabla_g\chi^{a}_\delta(x), \;\nabla_g\left[\hat{\delta}_{a,\lambda}(x)-\frac{c_0}{\sqrt{\lambda}d_g(a,x)}+\frac{c_0}{\sqrt{\lambda}d_g(a,x)}-\frac{\bar G_a(x)}{\sqrt{\lambda}}\right]\right>,\;\;\;\;x\in M.
\end{equation}
\vspace{4pt}
 
 \noindent
In the next step, using \eqref{bubble-M3} and \eqref{estgg}, we have the following:
\begin{equation}\label{diff3}
\left|\nabla_g\left[\hat{\delta}_{a,\lambda}(x)-\frac{c_0}{\sqrt{\lambda}d_g(a,x)}\right]\right|\le \frac{C}{\sqrt{\lambda}d_g(a,x)},\;\;\;\;x\in M, 
\end{equation}
and 
\begin{equation}\label{diff4}
\left|\nabla_g\left[\frac{c_0}{\sqrt{\lambda}d_g(a,x)}-\frac{\bar G_a(x)}{\sqrt{\lambda}}\right]\right|\le \frac{C}{\sqrt{\lambda}d_g(a,x)},\;\;\;\;x\in M.
\end{equation}
\vspace{4pt}
 
 \noindent
On the other hand, using \eqref{chi}, we derive 
\begin{equation}\label{fgchi}
\left|\nabla_g \chi^{a}_\delta(x)\right|\le \frac{C}{\delta} \large{1}_{\{y\in M:\;\delta\le d_g(a,y)\le 2\delta\}}(x),\;\;\;\;x\in M.
\end{equation}
Hence, combining \eqref{I_2} and \eqref{fgchi}, we get
\begin{equation}\label{estI2}
|I_2|\le \frac{C}{\delta^2\sqrt{\lambda}} \large{1}_{\{y\in M:\;\delta\le d_g(a,y)\le 2\delta\}}(x),\;\;\;\;x\in M.
\end{equation}
\vspace{6pt}
 
 \noindent
For \;$I_3,$\; using \eqref{chi} and \eqref{chid}, we obtain 
\begin{equation}\label{estI3}
    |I_3|\le C\hat{\delta}_{a,\lambda}(x) \large{1}_{\{y\in M:\;d_g(a,y)\le 2\delta\}}(x),\;\;\;\;x\in M.
\end{equation}
\vspace{4pt}
 
 \noindent
To deal with \;$I_4,$\; we first estimate \;$\Delta_g\hat{\delta}_{a,\lambda}$\; on \;$\psi_{a}(B^{+}_{2\delta}(0))$\; (see \eqref{psi}). For this, identifying \;$\psi_{a}(x)$\; with \;$x\in B^{+}_{2\delta}(0),$\; we have 
\begin{equation}\label{estI42}
\begin{split}
    \Delta_g\hat{\delta}_{a,\lambda}(x)&=\frac{1}{\sqrt{|g(x)|}}\partial_{i}\left[g^{ij}(x)\sqrt{|g(x)|}\partial_{j}\delta_{0,\lambda}(x)\right]+\frac{1}{\sqrt{|g(x)|}}\partial_{n}\left[\sqrt{|g(x)|}\partial_{n}\delta_{0,\lambda}(x)\right]\\
     &=\partial_{i}\left[g^{ij}(x)\frac{x_j}{r}\right]\partial_{r}\delta_{0,\lambda}(x)+g^{ij}(x)\frac{x_j}{r}\partial_{i}\partial_{r}\delta_{0,\lambda}(x)+\frac{1}{\sqrt{|g(x)|}}\partial_{i}\sqrt{|g(x)|}\frac{x_j}{r}g^{ij}(x)\partial_{r}\delta_{0,\lambda}(x)\\
     &+\frac{1}{\sqrt{|g(x)|}}\partial_{n}(\sqrt{|g(x)|})\partial_{n}\delta_{0,\lambda}(x)+\partial^{2}_{nn}\delta_{0,\lambda}(x)\\
     &=A+\Delta\delta_{0,\lambda}(x),\;\;\;\;\;\text{for}\;\;\;\;x\in B^+_{2\delta}(0),
    \end{split}
\end{equation}
where
\begin{equation}
\begin{split}
        A&=\left(g^{ij}(x)\frac{x_i x_j}{r^2}-1\right)\partial^{2}_{r}\delta_{0,\lambda}(x)+\left(g^{ij}(x)\frac{x_j}{r}\partial_{i}\log(\sqrt{|g(x)|})+\partial_{i}(g^{ij}(x)\frac{x_j}{r})-(\frac{n-1}{r})\right)\partial_{r}\delta_{0,\lambda}(x)\\
        &+\partial_{n}\log(\sqrt{|g(x)|})\partial_{i}\delta_{0,\lambda}(x),
        \end{split}
\end{equation}
and
\begin{equation}
        \D\delta_{0,\lambda}(x)=\partial^{2}_{r}\delta_{0,\lambda}(x)+(\frac{n-1}{r})\partial_{r}\delta_{0,\lambda}(x)+\partial^{2}_{n}\delta_{0,\lambda}(x).
   \end{equation}
   \vspace{4pt}
 
 \noindent
In \eqref{estI42}, \;$n=3,$\;
$$\partial_{i}=\frac{\partial}{\partial x_{i}}\;\;\;\;\;i=1,..,n,\;\;\;\;\text{and}\;\;\;\;\partial_{r}=\frac{\partial}{\partial r}.$$\;\\
 \vspace{4pt}
 
 \noindent
Since 
$$\Delta\delta_{0,\lambda}(x)=0,\;\;\;\;\text{in}\;\;\;\;B^+_{2\delta}(0)\cap\R^{3}_{+},$$\;then \eqref{estI42} implies
 \begin{equation}\label{estI43}
 \begin{split}
     \Delta_g\hat{\delta}_{a,\lambda}(x)&=\left(g^{ij}(x)\frac{x_i x_j}{r^2}-1\right)\partial^{2}_{r}\delta_{0,\lambda}(x)+\left(g^{ij}(x)\frac{x_j}{r}\partial_{i}\log(\sqrt{|g(x)|})+\partial_{i}(g^{ij}(x)\frac{x_j}{r})-(\frac{n-1}{r})\right)\partial_{r}\delta_{0,\lambda}(x)\\
     &+\partial_{n}\log(\sqrt{|g(x)|})\partial_{i}\delta_{0,\lambda}(x),\;\;\;\;x\in B^+_{2\delta}(0).
     \end{split}
\end{equation}
Combining \eqref{E_Metric} and \eqref{estI43}, we get
 \begin{equation}\label{estI44}
    |\Delta_g\hat{\delta}_{a,\lambda}(x)|\le C \left[|x|^{2}\partial^{2}_{r}\delta_{0,\lambda}(x)+|x|\left(\partial_{r}\delta_{0,\lambda}(x)+\partial_{n}\delta_{0,\lambda}(x)\right)\right],\;\;\;\;x\in B^+_{2\delta}(0). 
\end{equation}
Using Lemma \ref{Appendix} in \eqref{estI44}, we obtain the following estimate for \;$I_4$\; 
\begin{equation}\label{estI4}
    |I_4|\le C \left(\frac{\lambda}{1+\lambda^{2}d^{2}_g(a,x)}\right)^{\frac{1}{2}}\large{1}_{\{y\in M:\;d_g(a,y)\le 2\delta\}}(x),\;\;\;\;x\in M.
\end{equation}
Hence, the result for \eqref{C_01} follows from  \eqref{estI1}, \eqref{estI2}, \eqref{estI3}, and \eqref{estI4}.\\
\vspace{6pt}
 
\noindent
For the formula \eqref{C_02}, we first use \eqref{CC_0} to derive 
\begin{equation*}
\begin{split}
    -\frac{\partial u_{a,\lambda}(x)}{\partial n_{g}}-u^3_{a,\lambda}(x)&=\left(-\frac{\partial}{\partial n_{g}}\right)\left[\chi^{a}_\delta(x)\hat{\delta}_{a,\lambda}(x)+\left(1-\chi^{a}_\delta(x)\right)\frac{\bar{G}_a(x)}{\sqrt{\lambda}}\right]-u^3_{a,\lambda}(x)\\
    &=-\frac{\partial}{\partial n_{g}}\left[\chi^{a}_\delta(x)\hat{\delta}_{a,\lambda}(x)\right]-\frac{\partial}{\partial n_{g}}\left[\left(1-\chi^{a}_\delta(x)\right)\frac{\bar{G}_a(x)}{\sqrt{\lambda}}\right]-u^3_{a,\lambda}(x),\;\;\;\;x\in \partial M.
    \end{split}
\end{equation*}
To continue, we write 
\begin{equation*}
 -\frac{\partial}{\partial n_{g}}\left[\left(1-\chi^{a}_\delta(x)\right)\frac{\bar{G}_a(x)}{\sqrt{\lambda}}\right]=-\frac{\partial}{\partial n_{g}}\frac{\bar{G}_a(x)}{\sqrt{\lambda}}\left(1-\chi^{a}_\delta(x)\right)+\frac{\partial}{\partial n_{g}}\chi^{a}_\delta(x)\frac{\bar{G}_a(x)}{\sqrt{\lambda}},\;\;\;\;x\in \partial M.
\end{equation*}
Next, using the definition\;$\chi^{a}_\delta$\; (see \eqref{chi1}), the symmetry of \;$\chi^{a}_\delta$\; in \;$\psi_{a}(B^{+}_{2\delta}(0))$\; after passing to Euclidean coordinates, and \;$\frac{\partial G(a,x)}{\partial n_{g}}=0$\; for \;$x\in \partial M,$\;\;$x\neq a,$\; we have
\begin{equation*}
   -\frac{\partial}{\partial n_{g}}\left[\left(1-\chi^{a}_\delta(x)\right)\frac{\bar{G}_a(x)}{\sqrt{\lambda}}\right]=0,\;\;\;\;\text{for}\;\;\;\;x\in \partial M,\;\;\;\;\text{and}\;\;\;\;x\neq a.  
\end{equation*}
Hence for \;$x\in \partial M$\; and \;$x\neq a,$\; we have 
\begin{equation*}
 -\frac{\partial u_{a,\lambda}(x)}{\partial n_{g}}-u^3_{a,\lambda}(x)=-\frac{\partial}{\partial n_{g}}\left[\chi^{a}_\delta(x)\hat{\delta}_{a,\lambda}(x)\right]-u^3_{a,\lambda}(x).  
\end{equation*}
Using again the definition of \;$\chi^{a}_\delta$\; and the symmetry of \;$\chi^{a}_\delta$\; in  \;$\psi_{a}(B^{+}_{2\delta}(0))$\; as before, we obtain 
\begin{equation}\label{C_022}
 -\frac{\partial u_{a,\lambda}(x)}{\partial n_{g}}-u^3_{a,\lambda}(x)=-\chi^{a}_\delta(x)\frac{\partial \hat{\delta}_{a,\lambda}(x)}{\partial n_{g}}-u^3_{a,\lambda}(x),\;\;\;\;x\in \partial M,\;\;\;x\neq a.
\end{equation}
Clearly, \eqref{C_022} is true for \;$x=a.$\; On other hand, since identifying \;$x\in B^{+}_{2\delta}(0)$\; with \;$\psi_{a}(x),$\; we have  
\begin{equation*}
-\frac{\partial \hat{\delta}_{a,\lambda}(x)}{\partial n_{g}}=-\frac{\partial \delta_{0,\lambda}(x)}{\partial x_{3}}=\delta^{3}_{0,\lambda}(x),\;\;\;\;x\in \partial B^{+}_{2\delta}(0),
\end{equation*}
then  
\begin{equation*}
     -\frac{\partial u_{a,\lambda}(x)}{\partial n_{g}}-u^3_{a,\lambda}(x)=0,\;\;\;\;\forall\;\;x\in \partial B^{+}_{\delta}(0).
\end{equation*}
Therefore, \eqref{C_022} implies 
\begin{equation*}
   \left|-\frac{\partial u_{a,\lambda}(x)}{\partial n_{g}}-u^3_{a,\lambda}(x)\right|\le C \left(\frac{\lambda}{1+\lambda^{2}d^{2}_g(a,x)}\right)^{\frac{3}{2}}\large{1}_{\{y\in \partial M:\;d_g(a,y)\geq \delta\}}(x),\;\;\;\;x\in \partial M. 
\end{equation*}
\end{pf}
%
%
%
%
\section{PS-sequences and Deformation Lemma}\label{PSS}
In this Section, we discuss the asymptotic behavior of Palais-Smale (PS)-sequences for \;$J_q.$\; We also introduce the neighborhoods of potential critical points at infinity  of $J_q$ and their associated selection maps. As in other applications of the Barycenter technique of Bahri-Coron\cite{bc} (see \cite{aldawood}, \cite{jfe2}, and \cite{moa}), we also recall the associated Deformation Lemma.
\vspace{4pt}

\noindent
By some arguments which are classical by now see for example (\cite{sma2}, \cite{sma3}, \cite{fcm2}, and \cite{sma}), we have the following the profile decomposition for (PS)-sequences of \;$J_q.$\;
\begin{lem}\label{psseq}
Suppose that \;$(u_k)\subset H^{1}_{+}(M)$ is a (PS)-sequences for \;$J_q$, that is \;$\n J_q(u_k) \rightarrow 0$\; and \;$J_q(u_k)\rightarrow c$ \;up to a subsequence, and $\oint_{\partial M} u_k^4\;\;dS_g=c^{2}$\; for \;$k\in \N^{*},$\; then up to a subsequence, we have have there exists \;$u_{\infty}\geq 0,$\; an integer \;$p\geq 0,$\; a sequence of points \;${a_{i, k}}\in \partial M,\;\;i=1, \cdots, p,$\; and a sequence of positive numbers \;${\l_{i, k}},$\;\;$i=1, \cdots p,$\; such that\\
1)\\
$$
\left\{
\begin{split}
-\D_g u_{\infty}+qu_{\infty}&=0& \text{in}\;\;&M,\\
-\frac{\partial u_{\infty}}{\partial n_{g}}&=u^{3}_{\infty}&\text{on}\;\;&\partial M.
\end{split}
\right.
$$
2)\\
$$
||u_k-u_{\infty}-\sum_{i=1}^p u_{a_{i, k}, \l_{i, k}}||_q\longrightarrow 0,\;\; \text{as}\;\; k \longrightarrow \infty .
$$
3)\\
$$
J_q(u_k)^{2}\longrightarrow J_q(u_{\infty})^{2}+pS^{2},\;\; \text{as}\;\; k \longrightarrow \infty .
$$
4)\\\
For $i\neq j$, 
$$ 
\frac{\l_{i, k}}{\l_{j, k}}+\frac{\l_{j, k}}{\l_{i, k}}+\l_{i, k}\l_{j, k}G^{-2}(a_{i, k}, a_{j k})\longrightarrow +\infty,\;\; \text{as}\;\; k \longrightarrow \infty, 
$$
where \;$G$\; is as in \eqref{eqgreen}, and \;$\|\;\|_{q}$\; is as in \eqref{uqsqrt}.\\
\end{lem}
\vspace{4pt}

\noindent
To introduce the neighborhoods of potential critical points at infinity of \;$J_q$, we first fix
\begin{equation}\label{varepsilon0}
\varepsilon_0>0\;\;\; \text{and}\;\;\;\varepsilon_0 \simeq 0.
\end{equation}
 Furthermore, we choose
 \begin{equation}\label{nu0}
 \nu_0>1\;\;\;\text{and}\;\;\;\nu_0\simeq 1.
 \end{equation}
Then for \;$p\in \N^*$,\; and \;$0<\varepsilon\leq \varepsilon_0$, we define \;$V(p, \varepsilon)$\; the \;$(p, \varepsilon)$-neighborhood of potential critical points at infinity of \;$J_q$\; by
\begin{equation*}
\begin{split}
V(p, \varepsilon):=\{u\in H^{1}_{+}(M):&\;\;\exists\; a_1, \cdots, a_{p}\in \partial M,\;\;\alpha_1, \cdots, \alpha_{p}>0,\;\;\l_1, \cdots,\l_{p}>0,\;\;\l_i\geq \frac{1}{\varepsilon}\;\;\text{for}\;\;i=1\cdots, p,\\ 
&\Vert u-\sum_{i=1}^{p}\alpha_i u_{a_i, \l_i}\Vert_q\leq \varepsilon,\;\;\frac{\alpha_i}{\alpha_j}\leq \nu_0\;\;\text{and}\;\;\varepsilon_{i, j}\leq \varepsilon\;\;\text{for}\;\;i\neq j=1, \cdots, p\}.
\end{split}
\end{equation*}
\vspace{6pt}

\noindent
Concerning the sets \;$V(p, \varepsilon)$, for every \;$p\in \N^*$ \;there exists \;$0<\varepsilon_p\leq\varepsilon_0$\; such that for every \;$0<\varepsilon\leq \varepsilon_p$, we have
\begin{equation}\label{eq:mini}
\begin{cases}
\forall u\in V(p, \varepsilon)\;\; \text{the minimization problem}\;\;\min_{B_{\varepsilon}^{p}}\Vert u-\sum_{i=1}^{p}\alpha_iu_{a_i, \l_i}\Vert_q \\
\text{has a solution }\;(\bar \alpha, A, \bar \l)\in B_{\varepsilon}^{p} , \text{which is unique up to permutations,}
\end{cases}
\end{equation}
where \;$B^{p}_{\varepsilon}$\; is defined as 
\begin{equation*}
\begin{split}
B_{\varepsilon}^{p}:=\{(&\bar\alpha=(\alpha_1, \cdots, \alpha_p),\;\; A=(a_1, \cdots, a_p),\;\; \bar \l=\l_1, \cdots,\l_p))\in \R^{p}_+\times (\partial M)^p\times (0, +\infty)^{p}\\
&\frac{\alpha_i}{\alpha_j}\leq \nu_0\;\; \text{and}\;\; \varepsilon_{i, j}\leq \varepsilon,\;\;i\neq j=1, \cdots, p\}.
\end{split}
\end{equation*}
We define the selection map \;$s_p$\; via
\begin{equation*}
s_{p}: V(p, \varepsilon)\longrightarrow (\partial M)^p/\sigma_p\\
:
u\longrightarrow s_{p}(u)=A
\;\,\text{and} \,\;A\;\;\text{is given by}\;\,\eqref{eq:mini}.
\end{equation*}
To state the Deformation Lemma needed for the application of the algebraic topological argument of Bahri-Coron\cite{bc} for existence, we first  set
\begin{equation}\label{dfwp}
W_p:=\{u\in H^{1}_{+}(M)\;:\;J_q(u)\leq (p+1)^{\frac{1}{2}}\mathcal{S}\},
\end{equation}
 for \;$p\in \N,$\; where \;$\mathcal{S}$\; is as in \eqref{S}. 
\vspace{4pt}

\noindent
As in \cite{aldawood}, \cite{martndia2}, and \cite{nss}, we have Lemma \ref{psseq} implies the following Deformation Lemma.
\begin{lem}\label{deform}
Assuming that \;$J_q$ \;has no critical points, then for every \;$p\in \N^*$, up to taking \;$\varepsilon_p$\; given by \eqref{eq:mini} smaller, we have that for every\; $0<\varepsilon\leq \varepsilon_p$, the topological pair \;$(W_p,\; W_{p-1})$\; retracts by deformation onto \;$(W_{p-1}\cup A_p, \;W_{p-1})$\; with \;$V(p, \;\tilde \varepsilon)\subset A_p\subset V(p, \;\varepsilon)$\; where \;$0<\tilde \varepsilon<\frac{\varepsilon}{4}$\; is a very small positive real number and depends on \;$\varepsilon$.
\end{lem}
%
%
%
%
\section{Self-action estimates}\label{SAE}
In this Section, we establish some sharp self-action estimates that are essential for the use of the Bahri-Coron\cite{bc}'s Barycenter technique for existence. We start with the numerator of the functional \;$J_q.$\;
\begin{lem}\label{num}
Assuming that \;$\theta>0$\; is small, then there exists \;$C>0$\; such that \;$\forall$\;\;$a \in \partial M,$\;\;$\forall$\;\;$0<2\delta<\delta_0$\; and \;$\forall$\;\;$0<\frac{1}{\l}\le\theta\delta,$\; we have \\
\[\int_{M}\left(\left|\nabla_g u_{a,\lambda}\right|^{2} +q u^{2}_{a,\lambda}\right)\;dV_g \le \oint_{\partial M}u^{4}_{a,\lambda}\;dS_g+\frac{C}{\lambda}\left[1+\delta+\frac{1}{\lambda\delta^{2}}\right],\]
where \;$\delta_0$\; is as in \eqref{delta0}.
\end{lem}
\begin{pf}
Setting  
\begin{equation*}
    I=\int_{M}\left(\left|\nabla_g u_{a,\lambda}\right|^{2} +q u^{2}_{a,\lambda}\right)\;dV_g,
\end{equation*}
we have by Green's first identity
\begin{equation*}
    I=\oint_{\partial M}u^{4}_{a,\lambda}\;dS_g+\underbrace{\int_{M}\left(-\Delta_gu_{a,\lambda}+qu_{a,\lambda}\right)u_{a,\lambda}\;dV_g}_{\mbox{$I_1$}}+\underbrace{\oint_{\partial M}\left[-\frac{\partial u_{a,\lambda}}{\partial_{n_g}}-u^{3}_{a,\lambda}\right]u_{a,\lambda}\;dS_g}_{\mbox{$I_2$}}.
\end{equation*}
We are going to estimate \;$I_1$\; and \;$I_2$\; separately. Using Lemma \ref{c0estimate}, we have for \;$I_1$\;  
\begin{equation*}
\begin{split}
    |I_1|&\le \int_{M}\left|-\Delta_g u_{a,\lambda}+qu_{a,\lambda}\right|u_{a,\lambda}\;dV_g\\
    &\le \frac{C}{\delta^2\sqrt{\lambda}}\int_{M}u_{a,\lambda}\Large{1}_{\{x\in M:\;\delta\le d_g(a,x)\le 2\delta\}}\;dV_g\\
    &+C\int_{M}\hat{\delta}_{a,\lambda}u_{a,\lambda}\large{1}_{\{x\in M:\;d_g(a,x)\le 2\delta\}}\;dV_g \\
    &+C\int_{M}\left(\frac{\lambda}{1+\lambda^{2}d^{2}_g(a,x)}\right)^{\frac{1}{2}}u_{a,\lambda}\large{1}_{\{x\in M:\;d_g(a,x)\le 2\delta\}}\;dV_g.
    \end{split}
\end{equation*}
For the first term on the right-hand side of the formula above, we observe
\begin{equation}\label{selfI11}
\begin{split}
    \int_{M}u_{a,\lambda}\large{1}_{\{x\in M:\;\delta\le d_g(a,x)\le 2\delta\}}\;dV_g&\le C \int_{\{x\in M:\;\delta\le d_g(a,x)\le 2\delta\}}\left(\frac{\lambda}{1+\lambda^{2}d^{2}_g(a,x)}\right)^{\frac{1}{2}}\;dV_g\\
    &\le C\int_{\{x\in M:\;\delta\le d_g(a,x)\le 2\delta\}} \frac{1}{\sqrt{\lambda}d_g(a,x)}\;dV_g\\
    &\le \frac{C}{\sqrt{\lambda}}\int_{B^{+}_{2\delta}(0)\setminus B^{+}_{\delta}(0)}\frac{1}{|x|}\;dx\\
    &\le\frac{C}{\sqrt{\lambda}}\int_{\delta}^{2\delta}r\;dr\\
    & \le C \frac{\delta^2}{\sqrt{\lambda}}.
    \end{split}
\end{equation}
For the second term, we get
\begin{equation}\label{selfI12}
    \begin{split}
    \int_{M}\hat{\delta}_{a,\lambda}u_{a,\lambda}\large{1}_{\{x\in M:\;d_g(a,x)\le 2\delta\}}\;dV_g&\le C\int_{\{x\in M:\;d_g(a,x)\le2\delta\}}\left(\frac{\lambda}{1+\lambda^{2}d^{2}_g(a,x)}\right)\;dV_g\\
    &\le C\int_{\{x\in M:\;d_g(a,x)\le2\delta\}}\frac{1}{\lambda d^{2}_g(a,x)}\;dV_g\\
    &\le \frac{C}{\lambda}\int_{B^{+}_{2\delta}(0)}\frac{1}{|x|^{2}}\;dx\\
    &\le \frac{C}{\lambda}\int_{0}^{2\delta}1\;dr\\
    &\le C\frac{\delta}{\lambda}. 
    \end{split}
\end{equation}
Now, for the final term, we have
\begin{equation}\label{selfI13}
\begin{split}
    \int_{M}\left(\frac{\lambda}{1+\lambda^{2}d^{2}_g(a,x)}\right)^{\frac{1}{2}}u_{a,\lambda}\large{1}_{\{x\in M:\;d_g(a,x)\le 2\delta\}}\;dV_g&\le C\int_{\{x\in M:\;d_g(a,x)\le2\delta\}}\left(\frac{\lambda}{1+\lambda^{2}d^{2}_g(a,x)}\right)\;dV_g\\
    &\le C \int_{\{x\in M:\;d_g(a,x)\le 2\delta\}} \frac{1}{\lambda d^{2}_g(a,x)}\;dV_g\\
    &\le \frac{C}{\lambda}\int_{B^{+}_{2\delta}(0)}\frac{1}{|x|^{2}}\;dx\\
    &\le \frac{C}{\lambda}\int_{0}^{2\delta}1\;dr\\
    &\le C\frac{\delta}{\lambda}.
    \end{split}
\end{equation}
Thus, combining \eqref{selfI11}-\eqref{selfI13}, we get
\begin{equation}\label{selfI1}
    |I_1|\le\frac{C}{\lambda}\left(1+\delta\right).
\end{equation}
Next, in the case of \;$I_2,$\; we have
\begin{equation*}
    \begin{split}
        |I_2|&\le \oint_{\partial M}\left[-\frac{\partial u_{a,\lambda}}{\partial n_{g}}-u^{3}_{a,\lambda}\right]u_{a,\lambda}\;dS_g\\
        &\le C\oint_{\partial M}\left(\frac{\lambda}{1+\lambda^{2}d^{2}_g(a,x)}\right)^{\frac{3}{2}}u_{a,\lambda}\Large{1}_{\{x\in \partial M:\;d_g(a,x)\geq \delta\}}\;dS_g.
    \end{split}
\end{equation*}
On the right-hand side of the above formula, we observe
\begin{equation}
    \begin{split}
        \oint_{\partial M}\left(\frac{\lambda}{1+\lambda^{2}d^{2}_g(a,x)}\right)^{\frac{3}{2}}u_{a,\lambda}\Large{1}_{\{x\in \partial M:\;d_g(a,x)\geq \delta\}}\;dS_g&\le C\oint_{\{x\in \partial M:\;d_g(a,x)\geq \delta\}}\left(\frac{\lambda}{1+\lambda^{2}d^{2}_g(a,x)}\right)^{2}\;dS_g\\
        &\le C\oint_{\{x\in \partial M:\;d_g(a,x)\geq \delta\}}\frac{1}{\lambda^{2}d^{4}_g(a,x)}\;dS_g\\
        &\le \frac{C}{\lambda^{2}}\left[\int_{
        \hat{B}_{\delta_{0}}(0)\setminus \hat{B}_{\delta}(0)}\frac{1}{|x|^{4}}\;dx+1\right]\\
        &\le \frac{C}{\lambda^{2}}\left[\int_{\R^{2}\setminus \hat{B}_{\delta}(0)}^{+\infty}r^{-3}\;dr+1\right]\\
        &\le \frac{C}{\lambda^{2}}\left[\int_{\delta}^{+\infty}r^{-3}\;dr+1\right]\\
        &\le C \frac{1}{\lambda^{2}\delta^{2}}.
    \end{split}
\end{equation}
Thus, we have
\begin{equation}\label{selfI2}
    |I_2|\le \frac{C}{\lambda}\left(\frac{1}{\lambda\delta^{2}}\right).
\end{equation}
As a result, by combining \eqref{selfI1} and \eqref{selfI2}, we obtain
\[\int_{M}\left(\left|\nabla_g u_{a,\lambda}\right|^{2} +q u^{2}_{a,\lambda}\right)\;dV_g\le \oint_{\partial M}u^{4}_{a,\lambda}\;dS_g+\frac{C}{\lambda}\left[1+\delta+\frac{1}{\lambda\delta^{2}}\right].\]
\end{pf}\\
\noindent
For the denominator of \;$J_q$, we have.
\begin{lem}\label{denom}
 Assuming that \;$\theta>0$\; is small, then there exists \;$C>0$\; such that \;$\forall$\;\;$a \in \partial M,$\;\;$\forall$\;\;$0<2\delta<\delta_0$\; and \;$\forall$\;\;$0<\frac{1}{\l}\le \theta\delta,$\; we have
\[\oint_{\partial M} u^4_{a,\lambda}\;dS_g=\int_{\R^2} \delta^4_{0,\lambda}\;dx+O\left(\frac{1}{\lambda^2\delta^2}\right),\]
where \;$\delta_0$\; is as in \eqref{delta0}.
\end{lem}

\begin{pf}
We have
\begin{equation}\label{denom1}
\begin{split}
\oint_{\partial M} u^4_{a,\lambda}\;dS_g&=\oint_{\{x\in \partial M:\;d_g(a,x)\le \delta\}} u^4_{a,\lambda}\;dS_g+\oint_{\{x\in \partial M:\;\delta<d_g(a,x)\le 2\delta\}} u^4_{a,\lambda}\;dS_g\\
&+\oint_{\{x\in \partial M:\;d_g(a,x)>2\delta\}} u^4_{a,\lambda}\;dS_g.
\end{split}
\end{equation}
We are going to estimate each terms of the right-hand side of the formula of \eqref{denom1}. To begin, we start with the first term, we have 
\begin{equation}\label{denom2}
\begin{split}
\oint_{\{x\in \partial M:\;d_g(a,x)\le \delta\}} u^4_{a,\lambda}\;dS_g&=\oint_{\{x\in \partial M:\;d_g(a,x)\le \delta\}} \hat{\delta}^4_{a,\lambda}\;dS_g\\
&=\int_{\hat{B}_{\delta}(0)}\delta^4_{0,\lambda}\;dx\\
&=\int_{\R^2}\delta^4_{0,\lambda}\;dx-\int_{\{x\in \R^{2}:\;|x|\geq\delta\}} \delta^4_{0,\lambda}\;dx\\
&=\int_{\R^2}\delta^4_{0,\lambda}\;dx+O\left(\frac{1}{\lambda^2\delta^2}\right).
\end{split}
\end{equation}
To get to the second term, we have
\begin{equation}\label{denom3}
\begin{split}
\oint_{\{x\in \partial M:\;\delta<d_g(a,x)\le 2\delta\}} u^4_{a,\lambda}\;dS_g&\le C\oint_{\{x\in \partial M:\;\delta<d_g(a,x)\le 2\delta\}}\left(\frac{\lambda}{1+\lambda^2d_g(a,x)^2}\right)^{2}\;dS_g\\
&\le C \oint_{\{x\in \partial M:\;\delta<d_g(a,x)\le 2\delta\}}\frac{1}{\lambda^{2}d^{4}_g(a,x)}\;dS_g\\ 
&\le \frac{C}{\lambda^{2}}\int_{\hat{B}_{2\delta}(0)\setminus\hat{B}_{\delta}(0)}\frac{1}{|x|^{4}}\;dx\\
&\le \frac{C}{\lambda^2}\int_{\delta}^{2\delta} r^{-3}\;dr\\&\le \frac{C}{\lambda^2\delta^2}.
\end{split}
\end{equation}
\vspace{4pt}

\noindent
Next, using \eqref{estg}, we estimate the final term as follows 
\begin{equation}\label{denom4}
\begin{split}
\oint_{\{x\in \partial M:\;d_g(a,x)>2\delta\}} u^4_{a,\lambda}\;dS_g&=\oint_{\{x\in \partial M:\;d_g(a,x)>2\delta\}} \left(\frac{c_0}{\sqrt{\lambda}}G_a(x)\right)^{4}dS_g\\&=\frac{C}{\lambda^2}\oint_{\{x\in \partial M:\;d_g(a,x)>2\delta\}} G^{4}_a(x)\;dS_g\\
&\le C \oint_{\{x\in \partial M:\;d_g(a,x)>2\delta\}} \frac{1}{\lambda^{2}d^4_g(a,x)}\;dS_g\\
&\le C\oint_{\{x\in \partial M:\;2\delta< d_g(a,x)\le \delta_{0}\}}\frac{1}{\lambda^{2}d^4_g(a,x)}\;dS_g\\
&+C\oint_{\{x\in \partial M:\;d_g(a,x)>2\delta_{0}\}}\frac{1}{\lambda^{2}d^4_g(a,x)}\;dS_g\\
&\le \frac{C}{\lambda^{2}}\left[\oint_{\{x\in \partial M:\;2\delta< d_g(a,x)\le \delta_{0}\}}\frac{1}{d^4_g(a,x)}\;dS_g+C\right]\\
&\le \frac{C}{\lambda^{2}}\left[\int_{\hat{B}_{\delta_{0}}(0)\setminus \hat{B}_{2\delta}(0)}\frac{1}{|x|^{4}}\;dx+C\right]\\
&\le \frac{C}{\lambda^{2}}\left[\int_{2\delta}^{+\infty} r^{-3}\;dr+C\right]\\
&\le\frac{C}{\lambda^2\delta^2}.
\end{split}
\end{equation}
Finally, combining, \eqref{denom1}-\eqref{denom4}, we have \\\
\[\oint_{\partial M} u^4_{a,\lambda}\;dS_g=\int_{\R^2} \delta^4_{0,\lambda}\;dx+O\left(\frac{1}{\lambda^2\delta^2}\right).\]
\end{pf}
\vspace{4pt}

\noindent
We establish now the $J_q$-energy estimate of \;$u_{a, \lambda}$\; needed for the application of the Barycenter technique
of Bahri-Coron\cite{bc}  for existence.
\vspace{4pt}
\begin{cor}\label{sharpenergy}
Assuming that \;$\theta>0$\; is small, then there exists \;$C>0$\; such that \;$\forall$\;\;$a\in \partial M,$\;\;$\forall$\;\;$0<2\delta<\delta_0$\; and \;$\forall$\;\;$0<\frac{1}{\l}\le \theta\delta,$\; we have

\[J_q(u_{a, \l})\leq \mathcal{S}\left( 1+ C\left[\frac{1}{\l}+\frac{\delta}{\l}+\frac{1}{\delta^2\l^2}\right]\right),\]
where \;$\delta_0$\; is as in \eqref{delta0}.
\end{cor}
\begin{pf}
It follows from the properties of \;$\d_{0, \l}$ (see \eqref{S3} \eqref{S}), Lemma \ref{num}, and Lemma \ref{denom}.
\end{pf}

%
%
%
%
\section{Interaction estimates}\label{IE}
Throughout this Section, we derive some sharp inter-action estimates required for the algebraic topological argument of Bahri-Coron\cite{bc} for existence . We start with the following technical inter-action estimates. We anticipate that it will provide the needed inter-action estimates between \;$e_{ij}$\; and \;$\epsilon_{ji},$\; see \eqref{epij} and \eqref{eij} for the definitions of \;$e_{ij}$\; and \;$\epsilon_{ji}.$
\begin{lem}\label{interact1}
Assuming that \;$\theta>0$ is small, then there exists \;$C>0$\; such that \;$\forall$\;\;$\ a_i, a_j\in \partial M,$\; \;$\forall$\;\;$0<2\delta<\delta_0,$\; and \;$\forall$\;\;$0<\frac{1}{\l_j},$\;\;$\frac{1}{\l_i}\leq \theta\delta,$\; we have
\begin{equation*}
    \begin{split}
     \int_{M}&\left|\left(-\Delta_g+q\right) u_{a_j,\lambda_j}\right|u_{a_i,\lambda_i}\;dV_g+\oint_{\partial M}\left|-\frac{\partial}{\partial_{n_g}}u_{a_j,\lambda_j}-u^{3}_{a_j,\lambda_j}\right|u_{a_i,\lambda_i}\;dS_g\\
     &\le C\left[\delta+\frac{1}{\lambda_j\delta^{2}}\right]\left(\frac{\lambda_{i}}{\lambda_{j}}+\lambda_i\lambda_j d^{2}_g(a_{i},a_{j})\right)^{-\frac{1}{2}},
    \end{split}
\end{equation*}
where \;$\delta_0$\; is as in \eqref{delta0}.
\end{lem}
\begin{pf}
Using Lemma \ref{c0estimate}, we get the following:
\begin{equation}\label{A_j}
\begin{split}
    \underbrace{\left|\left(-\Delta_g +q\right)u_{a_j,\lambda_j}(x)\right|}_{\mbox{$A_j$}}&\le C\left[\frac{1}{\delta^2\sqrt{\lambda_j}}\Large{1}_{\{y\in M:\;\delta\le d_g(a_j,y)\le 2\delta\}}(x)+\hat{\delta}_{a_j,\lambda_j}(x)\large{1}_{\{y\in M:\; d_g(a_j,y)\le 2\delta\}}(x)\right.\\
  &+\left.\left(\frac{\lambda_j}{1+\lambda^{2}_{j}d^{2}_g(a_j,x)}\right)^{\frac{1}{2}}\large{1}_{\{y\in M:\;d_g(a_j,y)\le 2\delta\}}(x)\right],\;\;x\in M,
  \end{split}
\end{equation}
and
\begin{equation}\label{B_j}
    \underbrace{\left|-\frac{\partial u_{a_j,\lambda_j}(x)}{\partial n_{g}}-u^{3}_{a_j,\lambda_j}(x)\right|}_{\mbox{$B_j$}}\le C\left[\left(\frac{\lambda_j}{1+\lambda^{2}_{j}d^{2}_g(a_j,x)}\right)^{\frac{3}{2}}\Large{1}_{\{y\in \partial M:\;d_g(a_j,y)\geq \delta\}}(x)\right],\;\;x\in \partial M.
\end{equation}
For \;$x\in\{y\in M:\;\;d_g(a_j,y)\leq 2\d\},$ we have 
$$\left(\frac{\lambda_j}{1+\lambda^{2}_{j}d^{2}_g(a_j,x)}\right)\geq\left(\frac{\lambda_j}{1+4\lambda^2_j \delta^2}\right)\geq\frac{1}{\lambda_j \delta^2}\left[1+O(\frac{1}{\lambda^2_j \delta^2})\right]\geq\frac{1}{2 \lambda_j\delta^2}.$$
This implies $$\frac{1}{\sqrt{\lambda_j} \delta}\leq C \left(\frac{\lambda_j}{1+\lambda^{2}_{j}d^{2}_g(a_j,x)}\right)^{\frac{1}{2}}.$$\\
Thus, using \eqref{A_j} and Lemma \ref{Appendix}, we get 
\begin{equation}\label{A_j1}
    A_j\le C\left(1+\frac{1}{\delta}\right) \left(\frac{\lambda_j}{1+\lambda^{2}_{j}d^{2}_g(a_j,x)}\right)^{\frac{1}{2}}\large{1}_{\{y\in M:\;d_g(a_j,y)\le 4\delta\}}(x),\;\;x\in M.
\end{equation}
For \;$B_j$\; given by \eqref{B_j}, we have
\begin{equation}\label{B_j1}
    B_j\le C \left(\frac{\lambda_j}{1+\lambda^{2}_{j}d^{2}_g(a_j,x)}\right)^{\frac{3}{2}}\large{1}_{\{y\in \partial M:\;d_g(a_j,y)\geq \delta\}}(x),\;\;x\in \partial M.
\end{equation}
Now, using \eqref{A_j1}, we obtain
\begin{equation}\label{intest1}
\int_{M} A_j\;u_{a_i},_{\lambda_i}\;dV_g\le C\left(1+\frac{1}{\delta}\right) \underbrace{\int_{\{x\in M:\;d_g(a_j,x)\le4\delta\}}\left(\frac{\lambda_j}{1+\lambda_j^2d^{2}_g(a_j,x)}\right)^{\frac{1}{2}}\left(\frac{\lambda_i}{1+\lambda_i^2d^{2}_g(a_i,x)}\right)^{\frac{1}{2}}\;dV_g}_{\mbox{$I_1$}}.
\end{equation}
Also, using \eqref{B_j1}, we have 
\begin{equation}
\oint_{\partial M} B_j\;u_{a_i},_{\lambda_i}\;dS_g\le  C\underbrace{\oint_{\{x\in \partial M:\;d_g(a_j,x)\geq\delta\}}\left(\frac{\lambda_j}{1+\lambda_j^2d^{2}_g(a_j,x)}\right)^{\frac{3}{2}}\left(\frac{\lambda_i}{1+\lambda_i^2d^{2}_g(a_i,x)}\right)^{\frac{1}{2}}\;dS_g}_{\mbox{$I_2$}}.
\end{equation}
In the next step, we are going to estimate \;$I_1$\; as follows\\
\begin{equation*}
    \begin{split}
        I_1&=\int_{\{x\in M:\;d_g(a_j,x)\le4\delta\}}\left(\frac{\lambda_j}{1+\lambda_j^2d^{2}_g(a_j,x)}\right)^{\frac{1}{2}}\left(\frac{\lambda_i}{1+\lambda_i^2d^{2}_g(a_i,x)}\right)^{\frac{1}{2}}\;dV_g\\
        &=\underbrace{\int_{\{x\in M:\;2d_g(a_i,x)\le\frac{1}{\lambda_j}+d_g(a_i,a_j)\}\cap\{x\in M:\;d_g(a_j,x)\le 4\delta\}}\left(\frac{\lambda_j}{1+\lambda_j^2d^{2}_g(a_j,x)}\right)^{\frac{1}{2}}\left(\frac{\lambda_i}{1+\lambda_i^2d^{2}_g(a_i,x)}\right)^{\frac{1}{2}}}_{\mbox{$I^{1}_1$}}\;dV_g\\
        &+\underbrace{\int_{\{x\in M:\;2d_g(a_i,x)>\frac{1}{\lambda_j}+d_g(a_i,a_j)\}\cap\{x\in M:\;d_g(a_j,x)\le 4\delta\}}\left(\frac{\lambda_j}{1+\lambda_j^2d^{2}_g(a_j,x)}\right)^{\frac{1}{2}}\left(\frac{\lambda_i}{1+\lambda_i^2d^{2}_g(a_i,x)}\right)^{\frac{1}{2}}}_{\mbox{$I^{2}_1$}}\;dV_g.
    \end{split}
\end{equation*}
Using the triangle inequality, we estimate \;$I^{1}_1$\; as follows
\begin{equation*}
\begin{split}
    I^{1}_1&\le C\int_{\{x\in M:\;d_g(a_i,x)\le 8\delta\}}\left(\frac{\lambda_j}{1+\lambda_j^2d^{2}_g(a_j,a_i)}\right)^{\frac{1}{2}}\left(\frac{\lambda_i}{1+\lambda_i^2d^{2}_g(a_i,x)}\right)^{\frac{1}{2}}\;dV_g\\
    &\le C\frac{\sqrt{\frac{\lambda_j}{\lambda_i}}}{\left(1+\lambda^2_jd^{2}_g(a_j,a_i)\right)^{\frac{1}{2}}}\int_{\{x\in M:\;d_g(a_i,x)\le 8\delta\}} \frac{1}{d_g(a_i,x)}\;dV_g\\
     &\le C\frac{\sqrt{\frac{\lambda_j}{\lambda_i}}}{\left(1+\lambda^2_jd^{2}_g(a_j,a_i)\right)^{\frac{1}{2}}} \int_{B^{+}_{8\delta}(0)} \frac{1}{|x|}\;dx.
\end{split}
\end{equation*}
So, for \;$I^{1}_1,$\; we have
\begin{equation}\label{I11A_j}
    I^{1}_1=O\left(\delta^{2}\left(\frac{\lambda_{i}}{\lambda_{j}}+\lambda_i\lambda_j d^{2}_g(a_{i},a_{j})\right)^{-\frac{1}{2}}\right).
\end{equation}
For \;$I^{2}_1,$\; we observe 
\begin{equation*}
\begin{split}
    I^{2}_1& \le C\int_{\{x\in M:\;d_g(a_i,x)\le 4\delta\}}\left(\frac{\lambda_j}{1+\lambda_j^2d^{2}_g(a_j,x}\right)^{\frac{1}{2}}\left(\frac{\lambda_i}{1+\lambda_i^2d^{2}_g(a_j,a_i)}\right)^{\frac{1}{2}}\left(\frac{\lambda_j}{\lambda_i}\right)\;dV_g\\
    &\le C\frac{\sqrt{\frac{\lambda_i}{\lambda_j}}\left(\frac{\lambda_j}{\lambda_i}\right)}{\left(1+\lambda^2_jd^{2}_g(a_j,a_i)\right)^{\frac{1}{2}}}\int_{\{x\in M:\;d_g(a_i,x)\le 4\delta\}} \frac{1}{d_g(a_j,x)}\;dV_g\\
    &\le C\frac{\sqrt{\frac{\lambda_i}{\lambda_j}}\left(\frac{\lambda_j}{\lambda_i}\right)}{\left(1+\lambda^2_jd^{2}_g(a_j,a_i)\right)^{\frac{1}{2}}} \int_{B^{+}_{4\delta}(0)} \frac{1}{|x|}\;dx.
\end{split}
\end{equation*}
Thus, we get for \;$I^{2}_1$\;
\begin{equation}\label{I12A_j}
    I^{2}_1=O\left(\delta^{2}\left(\frac{\lambda_{i}}{\lambda_{j}}+\lambda_i\lambda_j d^{2}_g(a_{i},a_{j})\right)^{-\frac{1}{2}}\right).
\end{equation}
Hence, using \eqref{I11A_j} and \eqref{I12A_j}, we obtain
\begin{equation}\label{I1A_j}
     I_1=O\left(\delta^{2}\left(\frac{\lambda_{i}}{\lambda_{j}}+\lambda_i\lambda_j d^{2}_g(a_{i},a_{j})\right)^{-\frac{1}{2}}\right).
\end{equation}
For the last step, we are going to estimate \;$I_2.$\; We first write the following for \;$I_2$\;
\begin{equation*}
    \begin{split}
        I_2&=\oint_{\{x\in\partial M:\;d_g(a_j,x)\geq\delta\}}\left(\frac{\lambda_j}{1+\lambda_j^2d^{2}_g(a_j,x)}\right)^{\frac{3}{2}}\left(\frac{\lambda_i}{1+\lambda_i^2d^{2}_g(a_i,x)}\right)^{\frac{1}{2}}\;dS_g\\
        &=\underbrace{\oint_{\{x\in\partial M:\;2d_g(a_i,x)\le\frac{1}{\lambda_j}+d_g(a_j,a_i\})\cap\{x\in\partial M:\;d_g(a_j,x)\geq\delta)\}}\left(\frac{\lambda_j}{1+\lambda_j^2d^{2}_g(a_j,x)}\right)^{\frac{3}{2}}\left(\frac{\lambda_i}{1+\lambda_i^2d^{2}_g(a_i,x)}\right)^{\frac{1}{2}}}_{\mbox{$I^{1}_2$}}\;dS_g\\
        &+\underbrace{\oint_{\{x\in\partial M:\;2d_g(a_i,x)>\frac{1}{\lambda_j}+d_g(a_j,a_i)\}\cap\{x\in\partial M:\;d_g(a_j,x)\geq\delta\}}\left(\frac{\lambda_j}{1+\lambda_j^2d^{2}_g(a_j,x)}\right)^{\frac{3}{2}}\left(\frac{\lambda_i}{1+\lambda_i^2d^{2}_g(a_i,x)}\right)^{\frac{1}{2}}}_{\mbox{$I^{2}_2$}}\;dS_g.
    \end{split}
\end{equation*}
\noindent
In order to estimate \;$I^{1}_2,$\; we set
\begin{equation*}
    \begin{split}
    \mathcal{A}&=\{x\in\partial M:2d_g(a_i,x)\le\frac{1}{\lambda_j}+d_g(a_j,a_i)\}\cap\{x\in\partial M:d_g(a_j,x)\geq\delta\},\\
    r_{ij}&=\frac{1}{2}\left(\frac{1}{\lambda_j}+d_g(a_j,a_i)\right),
    \end{split}
\end{equation*}
and observe
\begin{equation*}
\begin{split}
I^{1}_2&\le \frac{C}{ \lambda_j^{-\frac{3}{2}}}\oint_{\mathcal{A}}\left(\frac{1}{1+\lambda_j^2d^{2}_g(a_j,a_i)}\right)^{\frac{1}{2}}\left(\frac{\lambda_i}{1+\lambda_i^2d^{2}_g(a_i,x)}\right)^{\frac{1}{2}}\left(\frac{1}{1+\lambda_j^2d^{2}_g(a_j,a_i)}\right)\;dS_g\\
&\le C\left(\frac{1}{\sqrt{\lambda_i}}\right)\left(\frac{1}{\lambda_j^2 \delta^2}\right) \frac{\lambda_j^{\frac{3}{2}}}{\left(1+\lambda_j^2d^{2}_g(a_j,a_i)\right)^{\frac{1}{2}}}\oint_{\{x\in\partial M:\;2d_g(a_i,x)\le r_{ij}\}}\frac{1}{d_g(a_i,x)}\;dS_g\\
&\le C\left(\frac{1}{\sqrt{\lambda_i}}\right)\left(\frac{1}{\lambda_j^2 \delta^2}\right) \frac{\lambda_j^{\frac{3}{2}}}{\left(1+\lambda_j^2d^{2}_g(a_j,a_i)\right)^{\frac{1}{2}}}\int_{\hat{B}_{r_{ij}}(0)}\frac{1}{|x|}\;dx\\
&\le C \left(\sqrt{\frac{\lambda_j}{\lambda_i}}\right)\left(\frac{1}{\lambda_j\delta^2}\right) \frac{1}{\left(1+\lambda_j^2d^{2}_g(a_j,a_i)\right)^{\frac{1}{2}}}\left(\frac{1}{\lambda_j}+d_g(a_j,a_i)\right)\\
&\le C \left(\sqrt{\frac{\lambda_j}{\lambda_i}}\right)\left(\frac{1}{\lambda_j \delta^{2}}\right) \frac{1}{\left(1+\lambda_j^2d^{2}_g(a_j,a_i)\right)^{\frac{1}{2}}}.
\end{split}
\end{equation*}
This implies 
\begin{equation}\label{I12B_j}
I^{1}_2=O\left(\frac{1}{\lambda_j\delta^{2}}\left(\frac{\lambda_{i}}{\lambda_{j}}+\lambda_i\lambda_j d^{2}_g(a_{i},a_{j})\right)^{-\frac{1}{2}}\right).
\end{equation}
For \;$I^{2}_2,$\; we derive
\begin{equation*}
\begin{split}
I^{2}_2&\le C\frac{\lambda_j}{\sqrt{\lambda_i}}\oint_{\{x\in \partial M:\;d_g(a_{j},x)\geq\delta\}}\left(\frac{1}{1+\lambda_j^2d^{2}_g(a_{i},a_{j})}\right)^{\frac{1}{2}}\left(\frac{\lambda_j}{1+\lambda_j^2d^{2}_g(a_{j},x)}\right)^{\frac{3}{2}}\;dS_g\\
&\le C\frac{\lambda_j}{\sqrt{\lambda_i}}\left(\frac{1}{1+\lambda_j^2d^{2}_g(a_{i},a_{j})}\right)^{\frac{1}{2}}\frac{1}{\lambda_j^{\frac{3}{2}}}\oint_{\{x\in \partial M:\;d_g(a_{j},x)\geq\delta\}}\frac{1}{d^{3}_g(a_{j},x)}\;dS_g\\
&\le C\frac{\lambda_j}{\sqrt{\lambda_i}}\left(\frac{1}{1+\lambda_j^2d^{2}_g(a_{i},a_{j})}\right)^{\frac{1}{2}}\frac{1}{\lambda_j^{\frac{3}{2}}}\oint_{\{x\in \partial M:\;\delta \le d_g(a_{j},x)\le \delta_{0}\}}\frac{1}{d^{3}_g(a_{j},x)}\;dS_g\\
&+ C\frac{\lambda_j}{\sqrt{\lambda_i}}\left(\frac{1}{1+\lambda_j^2d^{2}_g(a_{i},a_{j})}\right)^{\frac{1}{2}}\frac{1}{\lambda_j^{\frac{3}{2}}}\oint_{\{x\in \partial M:\;d_g(a_{j},x)\geq \delta_{0}\}}\frac{1}{d^{3}_g(a_{j},x)}\;dS_g\\
&\le C\frac{\lambda_j}{\sqrt{\lambda_i}}\left(\frac{1}{1+\lambda_j^2d^{2}_g(a_{i},a_{j})}\right)^{\frac{1}{2}}\frac{1}{\lambda_j^{\frac{3}{2}}}\left[\int_{\hat{B}_{\delta_{0}}(0)\setminus \hat{B}_{\delta}(0)} \frac{1}{|x|^{3}}\;dx+C\right]\\
&\le C\left(\sqrt{\frac{\lambda_j}{\lambda_i}}\right)\left(\frac{1}{\lambda_j \delta}\right) \frac{1}{\left(1+\lambda_j^2d^{2}_g(a_{i},a_{j})\right)^{\frac{1}{2}}}.
\end{split}
\end{equation*}
Thus, we get for \;$I^{2}_2$\;
\begin{equation}\label{I22B_j}
I^{2}_2=O\left(\frac{1}{\lambda_j\delta}\left(\frac{\lambda_{i}}{\lambda_{j}}+\lambda_i\lambda_j d^{2}_g(a_{i},a_{j})\right)^{-\frac{1}{2}}\right).
\end{equation}\\
Hence, combining \eqref{I12B_j} and \eqref{I22B_j}, we have
\begin{equation}\label{I2B_j}
I_2=O\left(\frac{1}{\lambda_j\delta^{2}}\left(\frac{\lambda_{i}}{\lambda_{j}}+\lambda_i\lambda_j d^{2}_g(a_{i},a_{j})\right)^{-\frac{1}{2}}\right).
\end{equation}\\
Therefore, using \eqref{A_j1}, \eqref{B_j1}, \eqref{I1A_j}, and \eqref{I2B_j}, we obtain\\
\begin{equation*}
    \begin{split}
     \int_{M}&\left|\left(-\Delta_g+q\right) u_{a_j,\lambda_j}\right|u_{a_i,\lambda_i}\;dV_g+\oint_{\partial M}\left|-\frac{\partial u_{a_j,\lambda_j}}{\partial_{n_g}}-u^{3}_{a_j,\lambda_j}\right|u_{a_i,\lambda_i}\;dS_g\\
     &\le C\left[\delta+\frac{1}{\lambda_j\delta^{2}}\right]\left(\frac{\lambda_{i}}{\lambda_{j}}+\lambda_i\lambda_j d^{2}_g(a_{i},a_{j})\right)^{-\frac{1}{2}}.   
    \end{split}
\end{equation*}
Hence, the proof of Lemma \ref{interact1} is complete.
\end{pf}\\
\vspace{2pt}

\noindent
 Clearly Lemma \ref{interact1} implies the following sharp interaction-estimate relating \;$e_{ij},$\;\;$\epsilon_{ij},$\; and \;$\varepsilon_{ij}$\; (for their
definitions, see \eqref{varepij}-\eqref{eij}).
\begin{cor}\label{interact2}
Assuming that \;$\theta>0$\; is small and \;$\mu_0>0$\; is small, then \;$\forall$\;\;$a_i, a_j\in \partial M$,\;\;$\forall 0<2\delta<\delta_0,$\; and \;$\forall$\;\;$0<\frac{1}{\l_j},$\;\;$\frac{1}{\l_i}\leq \theta\delta$\; such that \;$\varepsilon_{ij}\leq \mu_0,$\; we have 

\[e_{ij}=\epsilon_{ij}+O\left(\delta+\frac{1}{\lambda_i\delta^{2}}\right)\varepsilon_{ij},\]
where \;$\delta_0$\; is as in \eqref{delta0}.
\end{cor}
\vspace{6pt}

\noindent
 The following lemma gives a refined inter-action estimate relating \;$\epsilon_{ji}$\; and \;$\varepsilon_{ij}$.
 
\begin{lem}\label{interact3}
Assuming that \;$\theta>0$\; is small and \;$\mu_0>0$\; is small, then \;$\forall$\;\;$ a_i, a_j\in \partial M$,\;\;$\forall 0<2\delta<\delta_0,$\; and \;$\forall$\;\;$0<\frac{1}{\l_j}\le\frac{1}{\l_i}\leq \theta\delta$\; such that \;$\varepsilon_{ij}\leq \mu_0,$\; we have 
\[\epsilon_{ji}=c_0^4c_1\varepsilon_{ij}\left[\left(1+O\left(\delta+\frac{1}{\lambda^2_i\delta^2}\right)\right)\left(1+o_{\varepsilon_{ij}}(1)+O(\varepsilon_{ij} \delta^{-1})\right)+O\left(\varepsilon_{ij}\frac{1}{\delta^4}\right)\right],\]
where \;$c_0$\; is as in \eqref{bubble-M3} and \;$c_1$\; is as in \eqref{c3} .
\end{lem}

\begin{pf}
By using \eqref{uald}, we have
\[u_{a_i,\lambda_i}(x)=\chi^{a_i}_\delta(x)\hat{\delta}_{a_i,\lambda_i}(x)+(1-\chi^{a_i}_\delta(x))\frac{c_0}{\sqrt{\lambda}}G_{a_i}(x),\;\;\;x\in \partial M,\]
with \;$G_{a_i}(x)=G(a_i,x).$\; On the other hand, by using the definition of \;$\hat{\delta}_{a, \l}$\; (see \eqref{deltahat}), we have 
\[\chi^{a_i}_\delta(x)\hat{\delta}_{a_i,\lambda_i}(x)=c_0\chi^{a_i}_\delta(x)\left[\frac{\lambda_i}{1+\lambda^2_i G^{-2}_{a_i}(x)\frac{d^{2}_g(x,a_i)}{G^{-2}_{a_i}(x)}}\right]^{\frac{1}{2}},\;\;\;x\in \partial M.\]\\
\vspace{6pt}

\noindent
For \;$x\in\{y\in \partial M:\;d_g(a_i,y)\le 2\delta\},$\; the quantity \;$\left[1+\lambda^2_i G^{-2}_{a_i}(x)\frac{d^{2}_g(x,a_i)}{G^{-2}_{a_i}(x)}\right]$\; can be estimated by \\
\begin{equation*}
\begin{split}
1+\lambda^2_i G^2_{a_i}(x)\frac{d^2_g(x,a_i)}{G^{-2}_{a_i}(x)}&=1+\lambda^2_i G^{-2}_{a_i}(x)\left(1+O(\delta)\right)\\
&=1+\lambda^2_i G^{-2}_{a_i}(x)+O\left(\lambda^2_i \delta G^{-2}_{a_i}(x)\right)\\
&=\left(1+\lambda^2_i G^{-2}_{a_i}(x)\right)\left[1+O\left(\frac{\lambda^2_i \delta G^{-2}_{a_i}(x)}{1+\lambda^2_i G^{-2}_{a_i}(x)}\right)\right]\\
&=\left(1+\lambda^2_i G^{-2}_{a_i}(x)\right)\left[1+O(\delta)\right].
\end{split}
\end{equation*}\\
Hence, we have
\begin{equation}\label{part1}
\begin{split}
\chi^{a_i}_\delta(x) \delta_{a_i},_{\lambda_i}(x)&=c_0\chi^{a_i}_\delta(x)\left[\frac{\lambda_i}{\left(1+\lambda^2_i G^{-2}_{a_i}(x)\right)\left[1+O(\delta)\right]}\right]^{\frac{1}{2}}\\
&=c_0 \chi^{a_i}_\delta(x)\left[1+O(\delta)\right]\left[\frac{\lambda_i}{1+\lambda^2_i G^{-2}_{a_i}(x)}\right]^{\frac{1}{2}},\;\;\;x\in \partial M.
\end{split}
\end{equation}
Furthermore, we have
\[c_0(1-\chi^{a_i}_\delta(x))\left[\frac{\lambda_i}{1+\lambda^2_i G^{-2}_{a_i}(x)}\right]^{\frac{1}{2}}=(1-\chi^{a_i}_\delta(x))\frac{c_0}{\sqrt{\lambda_i}}G_{a_i}(x)\left[\frac{1}{1+\lambda^{-2}_i G^{2}_{a_i}(x)}\right]^{\frac{1}{2}},\;\;x\in \partial M.\]
Since on $\{x\in \partial M:\;d_g(x,a_i)\geq \d\}$, we have
\[ \frac{1}{1+\lambda^{-2}_iG^{2}_{a_i}(x)}=1+O\left(\frac{G^2_{a_i}(x)}{\lambda^2_i}\right)=1+O\left(\frac{1}{\lambda^2_i \delta^2}\right),\]
then we get
\[c_0(1-\chi^{a_i}_\delta(x))\left[\frac{\lambda_i}{1+\lambda^2_i G^{-2}_{a_i}(x)}\right]^{\frac{1}{2}}=(1-\chi^{a_i}_\delta(x))\frac{c_0}{\sqrt{\lambda_i}}G_{a_i}(x)\left(1+O\left(\frac{1}{\lambda^2_i \delta^2}\right)\right),\;\;\;x\in \partial M.\]
This implies 
\begin{equation}\label{part2}
(1-\chi^{a_i}_\delta(x))\frac{c_0}{\sqrt{\lambda}}G_{a_i}(x)=c_0(1-\chi^{a_i}_\delta(x))\left[\frac{\lambda}{1+\lambda^2_i G^{-2}_{a_i}(x)}\right]^{\frac{1}{2}}\left(1+O\left(\frac{1}{\lambda^2_i \delta^2}\right)\right),\;\;\;x\in \partial M.
\end{equation}
Thus, combining \eqref{part1} and \eqref{part2}, we get
\begin{equation*}
    \begin{split}
       u_{a_i},_{\lambda_i}(x)=c_0\left[\left(1+O(\delta)\right)\chi^{a_i}_\delta+(1-\chi^{a_i}_\delta)\left(1+O\left(\frac{1}{\lambda^2_i \delta^2}\right)\right)\right]\left[\frac{\lambda}{1+\lambda^2_i G^{-2}_{a_i}(x)}\right]^{\frac{1}{2}},\;\;\;x\in \partial M. 
    \end{split}
\end{equation*}
Hence, we obtain
\begin{equation}\label{partf}
u_{a_i, \lambda_i}(x)=c_0\left[1+O(\delta)+O\left(\frac{1}{\lambda^2_i\delta^2}\right)\right]\left[\frac{\lambda}{1+\lambda^2_i G^{-2}_{a_i}(x)}\right]^{\frac{1}{2}},\;\;\;x\in \partial M.
\end{equation}\\
\vspace{4pt}
\noindent
Now, we are going to complete our task by using \eqref{partf}. We begin by writing the estimate in the following form:
\begin{equation}\label{J1+J2}
\epsilon_{ji}=\oint_{\partial M} u^3_{a_j,\lambda_j}u_{a_i,\lambda_i}\;dS_g=\underbrace{\oint_{\hat{B}(a_j,\delta)} u^3_{a_j,\lambda_j}u_{a_i,\lambda_i}\;dS_g}_{\mbox{$J_1$}}+\underbrace{\oint_{\partial M-\hat{B}(a_j,\delta)} u^3_{a_j,\lambda_j}u_{a_i,\lambda_i}\;dS_g}_{\mbox{$J_2$}}.
\end{equation}
We are going to estimate\;$J_1$\;and\;$J_2$\; separately. In the case of \;$J_2,$\; we have
\begin{equation*}
\begin{split}
\oint_{\partial M-\hat{B}(a_j,\delta)} u^3_{a_j,\lambda_j}u_{a_i,\lambda_i}\;dS_g&\le C\oint_{\partial M-\hat{B}(a_j,\delta)}\left(\frac{1}{\lambda_j}\right)^{\frac{3}{2}}\left(\frac{1}{\delta}\right)^3u_{a_i,\lambda_i}\;dS_g\\
&\le C \left(\frac{1}{\lambda_j}\right)^{\frac{3}{2}}\left(\frac{1}{\delta}\right)^3\oint_{\partial M-\hat{B}(a_j,\delta)}u_{a_i,\lambda_i}\;dS_g\\
&\le C\left(\frac{1}{\lambda_j\d^2}\right)^{\frac{3}{2}}\oint_{\partial M-(\hat{B}(a_j,\delta)\cup \hat{B}(a_i,\delta))}u_{a_i,\lambda_i}\;dS_g\\
&+C\left(\frac{1}{\lambda_j\d^2}\right)^{\frac{3}{2}}\oint_{(\partial M-\hat{B}(a_j,\delta))\cap \hat{B}(a_i,\delta)}u_{a_i,\lambda_i}\;dS_g\\
&\le C\left(\frac{1}{\lambda_j}\right)^{\frac{3}{2}}\left(\frac{1}{\delta}\right)^4\frac{1}{\sqrt{\lambda_i}}+C\left(\frac{1}{\lambda_j}\right)^{\frac{3}{2}}\left(\frac{1}{\delta}\right)^3\oint_{\hat{B}(a_i,\delta)}\left(\frac{\lambda_i}{1+\lambda^2_id^{2}_g(a_i,x)}\right)^{\frac{1}{2}}\;dS_g\\
&\le C\left(\frac{1}{\lambda_j}\right)^{\frac{3}{2}}\left(\frac{1}{\delta}\right)^4\frac{1}{\sqrt{\lambda_i}}+C\delta\left(\frac{1}{\lambda_j}\right)^{\frac{3}{2}}\left(\frac{1}{\delta}\right)^3\frac{1}{\sqrt{\lambda_i}}\\
&\le C\left(\frac{1}{\lambda_j}\right)^{\frac{3}{2}}\left(\frac{1}{\delta}\right)^4\frac{1}{\sqrt{\lambda_i}}\left(1+\delta^2\right)\\
&\le \frac{C}{\lambda^{\frac{3}{2}}_j\sqrt{\lambda_i}\delta^4}.
\end{split}
\end{equation*}

\noindent
Thus, we  get for \;$J_2$\;
\begin{equation}\label{J_22}
J_2=O\left(\varepsilon^{2}_{ij}\frac{1}{\delta^4}\right).
\end{equation}

\noindent
In the next step, using \eqref{partf} for \;$J_1,$\; we have 
\begin{equation*}
\begin{split}
\oint_{\hat{B}(a_j,\delta)} u^3_{a_j,\lambda_j}u_{a_i,\lambda_i}\;dS_g&=c^4_0\oint_{\hat{B}(a_j,\delta)}\left(\frac{\lambda_j}{1+\lambda^2_jd^{2}_g(a_j,x)}\right)^{\frac{3}{2}}\left[1+O(\delta)+O\left(\frac{1}{\lambda^2_i\delta^2}\right)\right]\left[\frac{\lambda_i}{1+\lambda^2_i G^{-2}_{a_i}(x)}\right]^{\frac{1}{2}}\;dS_g\\
&=c^4_0\left[1+O(\delta)+O\left(\frac{1}{\lambda^2_i\delta^2}\right)\right]\oint_{\hat{B}(a_j,\delta)}\left(\frac{\lambda_j}{1+\lambda^2_jd^{2}_g(a_j,x)}\right)^{\frac{3}{2}}\left[\frac{\lambda_i}{1+\lambda^2_i G^{-2}_{a_i}(x)}\right]^{\frac{1}{2}}\;dS_g\\
&=c^4_0 \frac{1}{\sqrt{\lambda_j}}\left[1+O\left(\d+\frac{1}{\lambda^2_i\delta^2}\right)\right]\int_{\hat{B}_{\lambda_j\delta}(0)}\left(\frac{1}{1+|y|^2}\right)^{\frac{3}{2}}\left[\frac{\lambda_i}{1+\lambda^2_i G^{-2}_{a_i}\left(\psi_{a_j}\left(\frac{y}{\lambda_j}\right)\right)}\right]^{\frac{1}{2}}\\
&=c^4_0 \left[1+O\left(\d+\frac{1}{\lambda^2_i\delta^2}\right)\right]\int_{\hat{B}_{\lambda_j\delta}(0)}\left(\frac{1}{1+|y|^2}\right)^{\frac{3}{2}}\left[\frac{1}{\frac{\l_j}{\l_i}+\lambda_i \l_jG^{-2}_{a_i}\left(\psi_{a_j}\left(\frac{y}{\lambda_j}\right)\right)}\right]^{\frac{1}{2}}.
\end{split}
\end{equation*}
\vspace{6pt}

\noindent
Recalling that $\lambda_i\le\lambda_j$, then  for \ $\varepsilon_{ij}\sim 0 $, we have\\
1) Either $\varepsilon^{-2}_{ij}\sim \lambda_i\lambda_jG^{-2}_{a_i}(a_j).$
\vspace{4pt}

\noindent
2) or $\varepsilon^{-2}_{ij}\sim\frac{\lambda_j}{\lambda_i}$.
\vspace{6pt}

\noindent
In order to estimate \;$J_1,$\; we first define the following sets
\begin{equation*}
    \begin{split}
        A_1&=\left\{y\in \R^{2}:\;\left(\left|\frac{y}{\lambda_j}\right|\le \epsilon \:G^{-1}_{a_i}(a_j)\right)\cap \hat{B}_{\lambda_j\delta}(0)\right\},\\
        A_2&=\left\{y\in \R^{2}:\;\left(\left|\frac{y}{\lambda_j}\right|\le\epsilon\frac{1}{\lambda_i}\right)\cap \hat{B}_{\lambda_j\delta}(0)\right\},
    \end{split}
\end{equation*}
and
\[\mathcal{A}=A_1\cup A_2\;,\]
with \;$\epsilon>0$\; very small. Then  by Taylor expansion on \;$\mathcal{A}$, we have
\begin{equation*}
\begin{split}
\left[\frac{\lambda_j}{\lambda_i}+\lambda_i\lambda_jG^{-2}_{a_i}\left(\psi_{a_j}\left(\frac{y}{\lambda_j}\right)\right)\right]^{-\frac{1}{2}}&=\left[\frac{\lambda_j}{\lambda_i}+\lambda_i\lambda_jG^{-2}_{a_i}(a_j)\right]^{-\frac{1}{2}}\\
&+\left[\left(-\frac{1}{2}\nabla  G^{-2}_{a_i}\circ \psi_{a_j}(a_j)\lambda_i y\right)\right]\left[\frac{\lambda_j}{\lambda_i}+\lambda_i\lambda_jG^{-2}_{a_i}\left(\psi_{a_j}\left(\frac{y}{\lambda_j}\right)\right)\right]^{-\frac{3}{2}}\\
&+O\left[\left(\frac{\lambda_i}{\lambda_j}\right)|y|^2\right]\left[\frac{\lambda_j}{\lambda_i}+\lambda_i\lambda_jG^{-2}_{a_i}\left(\psi_{a_j}\left(\frac{y}{\lambda_j}\right)\right)\right]^{-\frac{3}{2}}.
\end{split}
\end{equation*}

\noindent
So, we write \;$J_1$\; such as  
\begin{equation}\label{J_11}
J_1=c_0^4\left[1+O(\delta)+O\left(\frac{1}{\lambda^2_i\delta^2}\right)\right]\left(\sum_{m=1}^{4} I_m\right),
\end{equation}
with
\begin{equation*}
\begin{split}
I_1&=\left[\frac{\lambda_j}{\lambda_i}+\lambda_i\lambda_jG^{-2}_{a_i}(a_j)\right]^{-\frac{1}{2}}\int_{\mathcal{A}}\left(\frac{1}{1+|y|^2}\right)^{\frac{3}{2}},\\
I_2&=\left[\frac{\lambda_j}{\lambda_i}+\lambda_i\lambda_jG^{-2}_{a_i}(a_j)\right]^{-\frac{3}{2}}\int_{\mathcal{A}}\left(\frac{1}{1+|y|^2}\right)^{\frac{3}{2}}\left[\nabla G^{-2}_{a_i}\circ \psi_{a_j}(a_j)\lambda_i y\right],\\
I_3&=\left[\frac{\lambda_j}{\lambda_i}+\lambda_i\lambda_jG^{-2}_{a_i}(a_j)\right]^{-\frac{3}{2}}\int_{\mathcal{A}}\left(\frac{1}{1+|y|^2}\right)^{\frac{3}{2}}O\left[\left(\frac{\lambda_i}{\lambda_j}\right)|y|^2\right],
\end{split}
\end{equation*}
and
\begin{equation*}
\begin{split}
I_4&=\int_{\hat{B}_{\lambda_j\delta}(0)-\mathcal{A}}\left(\frac{1}{1+|y|^2}\right)^{\frac{3}{2}}\left[\frac{\lambda_j}{\lambda_i}+\lambda_i\lambda_jG^{-2}_{a_i}\left(\psi_{a_j}\left(\frac{y}{\lambda_j}\right)\right)\right]^{-\frac{1}{2}}.
\end{split}
\end{equation*}

\noindent
Now, let us estimate \;$I_1.$\; We have 
\[I_1=\left[\frac{\lambda_j}{\lambda_i}+\lambda_i\lambda_jG^{-2}_{a_i}(a_j)\right]^{-\frac{1}{2}}\left[c_1+\int_{\R^2-\mathcal{A}}\left(\frac{1}{1+|y|^2}\right)^{\frac{3}{2}}\right],\]
where \;$c_1$\; is as in \eqref{c3}. On the other hand, we set
\begin{equation*}
    \begin{split}
     T_{ij}&=\lambda_j\epsilon G^{-1}_{a_i}(a_i,a_j),\\
     L_{ij}&=\epsilon\frac{\lambda_j}{\lambda_i},
    \end{split}
\end{equation*}
and have
\[\int_{\R^2-\mathcal{A}}\left(\frac{1}{1+|y|^2}\right)^{\frac{3}{2}} \le \int_{\R^2-\hat{B}_{\delta\lambda_j}(0)}\left(\frac{1}{1+|y|^2}\right)^{\frac{3}{2}}+\int_{\R^2-\hat{B}_{T_{ij}}(0)}\left(\frac{1}{1+|y|^2}\right)^{\frac{3}{2}}\]\\
\noindent
if \;$\varepsilon^{-2}_{ij}\sim \lambda_i\lambda_jG^{-2}_{a_i}(a_j),$\; and 
\begin{equation*}
\int_{\R^2-\mathcal{A}}\left(\frac{1}{1+|y|^2}\right)^{\frac{3}{2}}\le \int_{\R^2-\hat{B}_{\delta\lambda_j}(0)}\left(\frac{1}{1+|y|^2}\right)^{\frac{3}{2}}+\int_{\R^2-\hat{B}_{L_{ij}}(0)}\left(\frac{1}{1+|y|^2}\right)^{\frac{3}{2}}
\end{equation*}
\noindent
if \;$\varepsilon^{-2}_{ij}\sim\frac{\lambda_j}{\lambda_i}.$\; We have 
\begin{equation*}
\int_{\R^2-\hat{B}_{\delta\lambda_j}(0)}\left(\frac{1}{1+|y|^2}\right)^{\frac{3}{2}}=O\left(\frac{1}{\lambda_j\delta}\right).
\end{equation*}
\noindent
Moreover, if \;$\varepsilon^{-2}_{ij}\sim \lambda_i\lambda_jG^{-2}_{a_i}(a_j),$\; then
\begin{equation*}
    \begin{split}
\int_{\R^2-\hat{B}_{T_{ij}}(0)}\left(\frac{1}{1+|y|^2}\right)^{\frac{3}{2}}&=O\left(\frac{1}{\lambda_j G^{-1}_{a_i}(a_j)}\right)\\
&=O\left(\frac{1}{\sqrt{\lambda_j\lambda_i} G^{-1}_{a_i}(a_j)}\right)\\
&=O\left(\varepsilon_{ij}\right).
\end{split}
\end{equation*}

\noindent
Furthermore if \;$\varepsilon^{-2}_{ij}\sim\frac{\lambda_j}{\lambda_i},$ then
\begin{equation*}
\int_{\R^2-\hat{B}_{L_{ij}}(0)}\left(\frac{1}{1+|y|^2}\right)^{\frac{3}{2}}=O\left(\varepsilon_{ij}\right).
\end{equation*}
This implies 
\begin{equation*}
\int_{\R^2-\mathcal{A}}\left(\frac{1}{1+|y|^2}\right)^{\frac{3}{2}}=O\left(\varepsilon_{ij}+\frac{1}{\l_j\delta}\right)=O\left(\varepsilon_{ij}+\varepsilon_{ij}\frac{1}{\delta}\right)=O\left(\varepsilon_{ij}\frac{1}{\delta}\right).
\end{equation*}

\noindent
Thus, we get\
\begin{equation*}
    \begin{split}
I_1&=\left[\frac{\lambda_j}{\lambda_i}+\lambda_i\lambda_jG^{-2}_{a_i}(a_j)\right]^{-\frac{1}{2}}\left[c_1+O\left(\varepsilon_{ij}\frac{1}{\delta}\right)\right]\\
&=\varepsilon_{ij}\left(1+o_{\varepsilon_{ij}}(1)\right)\left[c_1+O\left(\varepsilon_{ij}\frac{1}{\delta}\right)\right]
\end{split}
\end{equation*}
\noindent
Hence, we obtain
\begin{equation}\label{I_11}
I_1=c_1\varepsilon_{ij}\left[1+o_{\epsilon_{ij}}(1)+O\left(\varepsilon_{ij}\frac{1}{\delta}\right)\right].
\end{equation}
By symmetry, we have
 \begin{equation}\label{I_22}
 I_2=0.
 \end{equation}
 
 \noindent
 Next, for \;$I_3,$\; we derive
\begin{equation*}
\begin{split}
\int_{\mathcal{A}}\frac{|y|^2}{\left(1+|y^2|\right)^{\frac{3}{2}}}
&\le \int_{\hat{B}_{T_{ij}}(0)}\frac{|y|^2}{\left(1+|y|^{2}\right)^{\frac{3}{2}}}+\int_{\hat{B}_{L_{ij}}(0)}\frac{|y|^2}{\left(1+|y|^2\right)^{\frac{3}{2}}}\\
&=O\left(\epsilon\lambda_jG^{-1}_{a_i}(a_j)+\epsilon \frac{\lambda_j}{\lambda_i}\right).
\end{split}
\end{equation*}
Thus, we have 
\begin{equation*}
\begin{split}
I_3&=\varepsilon^3_{ij}\left(\frac{\lambda_i}{\lambda_j}\right)\left(1+o_{\varepsilon_{ij}}(1)\right)\left[O\left(\epsilon\lambda_jG^{-1}_{a_i}(a_j)+\epsilon \frac{\lambda_j}{\lambda_i}\right)\right]\\
&=\varepsilon^3_{ij}\left(1+o_{\varepsilon_{ij}}(1)\right)\left[O\left(\sqrt{\lambda_i\lambda}_jG^{-1}_{a_i}(a_j)+\sqrt{\frac{\lambda_j}{\lambda_i}}\right)\right].
\end{split}
\end{equation*}

\noindent
Hence, we obtain
\begin{equation}\label{I_33}
I_3=O\left(\varepsilon^2_{ij}\right).
\end{equation}

\noindent
Finally, we estimate \;$I_4$\; as follows.\\\\
If \;$\varepsilon^{-2}_{ij}\sim\frac{\lambda_j}{\lambda_i},$\; then\\
\begin{equation}\label{part1i4}
I_4\le C\varepsilon_{ij}\int_{\hat{B}_{\lambda_j\delta}(0)-\mathcal{A}}\left(\frac{1}{1+|y|^2}\right)^{\frac{3}{2}}\le C\varepsilon_{ij}\left(\frac{\lambda_j}{\lambda_i}\right)^{-1}\le C \varepsilon^3_{ij}.
\end{equation}\\
If \;$\varepsilon^{-2}_{ij}\sim \lambda_i\lambda_jG^{-2}_{a_i}(a_j),$\; then we argue as follows.
In case \;$d_g(a_i,a_j)\geq2\delta,$\; since 
\[G_{a_i}\left(\psi_{a_j}\left(\frac{y}{\lambda_j}\right)\right)\le C\delta^{-1}\]
 for \;$y\in \hat{B}(0, \l_j\d),$\; then we have 
\begin{equation*}
\begin{split}
I_4&\le C\int_{\hat{B}_{\lambda_j\delta}(0)-\mathcal{A}}\left(\frac{1}{1+|y|^2}\right)^{\frac{3}{2}} \frac{1}{\sqrt{\lambda_i\lambda_j}\delta}\\
&\le\frac{C}{\sqrt{\lambda_i\lambda_j}\delta}\left(\frac{1}{\lambda_jG_{a_i}^{-1}(a_j)}\right)\\&\le \frac{C}{\sqrt{\lambda_i\lambda_j}G^{-1}_{a_i}(a_j)}\frac{G^{-1}_{a_i}(a_j)}{\sqrt{\lambda_i\lambda_j}G_{a_i}^{-1}(a_j)\d}\\&\le C\varepsilon_{ij} \varepsilon_{ij}\frac{1}{\delta} .
\end{split}
\end{equation*}

\noindent
Thus, when \; $d_g(a_i,a_j)\geq 2\d,$\; we have 
\begin{equation}\label{part2i4}
I_4=O\left(\varepsilon^2_{ij}\frac{1}{\delta}\right).
\end{equation}
In case \;$d_g(a_i,a_j)<2\delta,$\; we first observe that 
$$
\hat{B}_{\l_j\d}(0)\setminus \mathcal{A}\subset A_1\cup A_2
$$
with 
$$
A_1=\left\{y\in \R^{2}:\;\epsilon\lambda_jG^{-1}_{a_i}(a_j)\le|y|\le E\lambda_jG^{-1}_{a_i}(a_j)\right\}
$$
and
$$
A_2=\left\{y\in \R^{2}:\;E\lambda_jG^{-1}_{a_i}(a_j)\le|y|\le \lambda_j\delta\right\},
$$
where \;\;$0<\epsilon<E.$\;\\\\
\noindent
Thus, we have 
\begin{equation}\label{i4ine}
I_4\le I^1_4+I^2_4,
\end{equation}
with
\begin{equation*}
I_4^1=\int_{A_1} \left(\frac{1}{1+|y|^2}\right)^{\frac{3}{2}}\left[\frac{\lambda_j}{\lambda_i}+\lambda_i\lambda_jG^{-2}_{a_i}\left(\psi_{a_j}\left(\frac{y}{\lambda_j}\right)\right)\right]^{-\frac{1}{2}},
\end{equation*}
and
\begin{equation*}
I_4^2=\int_{A_2} \left(\frac{1}{1+|y|^2}\right)^{\frac{3}{2}}\left[\frac{\lambda_j}{\lambda_i}+\lambda_i\lambda_jG^{-2}_{a_i}\left(\psi_{a_j}\left(\frac{y}{\lambda_j}\right)\right)\right]^{-\frac{1}{2}}.
\end{equation*}
We estimate \;$I_4^1$\; as follows:
\begin{equation*}
\begin{split}
I_4^1&\le C\left[1+\lambda_j^2G^{-2}_{a_i}(a_j)\right]^{-\frac{3}{2}}\int_{\{y\in \R^{2}:\;|y|\le E\lambda_jG^{-1}_{a_i}(a_j)\}}\left[\frac{\lambda_j}{\lambda_i}+\lambda_i\lambda_jG^{-2}_{a_i}\left(\psi_{a_j}\left(\frac{y}{\lambda_j}\right)\right)\right]^{-\frac{1}{2}}\\
&\le C\left[1+\lambda_j^2G^{-2}_{a_i}(a_j)\right]^{-\frac{3}{2}}\left(\frac{\lambda_i}{\lambda_j}\right)^{\frac{1}{2}}\int_{\{y\in \R^{2}:\;|y|\le E\lambda_jG^{-1}_{a_i}(a_j)\}}\left[1+\lambda^2_iG^{-2}_{a_i}\left(\psi_{a_j}\left(\frac{y}{\lambda_j}\right)\right)\right]^{-\frac{1}{2}}\\
&\le C\left[1+\lambda_j^2G^{-2}_{a_i}(a_j)\right]^{-\frac{3}{2}}\left(\frac{\lambda_i}{\lambda_j}\right)^{\frac{1}{2}}\int_{\{y\in \R^{2}:\;|y|\le E\lambda_jG^{-1}_{a_i}(a_j)\}}\left[1+\lambda^2_i\left|\psi^{-1}_{a_i}\circ\psi_{a_j}\left(\frac{y}{\lambda_j}\right)\right|^2\right]^{-\frac{1}{2}}\\
&\le C\left[1+\lambda_j^2G^{-2}_{a_i}(a_j)\right]^{-\frac{3}{2}}\left(\frac{\lambda_i}{\lambda_j}\right)^{\frac{1}{2}}\left(\frac{\lambda_j}{\lambda_i}\right)^{2}\int_{\{z\in \R^{2}:\;|z|\le \bar E\lambda_iG^{-1}_{a_i}(a_j)\}}\left[\frac{1}{1+|z|^2}\right]^{\frac{1}{2}}\\
&\le C \left[\frac{\lambda_i}{\lambda_j}+\lambda_i\lambda_jG^{-2}_{a_i}(a_j)\right]^{-\frac{3}{2}}\left(\lambda_iG^{-1}_{a_i}(a_j)\right)
\\&
\le C \varepsilon^3_{ij}\left(\sqrt{\lambda_i\lambda_j} G^{-1}_{a_i}(a_j)\right),
\end{split}
\end{equation*}

\noindent
where \;$\bar E$\; is a positive constant. So we obtain
\begin{equation}\label{i14fin}
I^1_4=O\left(\varepsilon^2_{ij}\right).
\end{equation}
\noindent
For \;$I_4^2,$\; we have
\begin{equation*}
\begin{split}
I^2_4&=\int_{A_2}\left(\frac{1}{1+|y|^2}\right)^{\frac{3}{2}}\left[\frac{\lambda_j}{\lambda_i}+\lambda_i\lambda_jG^{-2}_{a_i}\left(\psi_{a_j}\left(\frac{y}{\lambda_j}\right)\right)\right]^{-\frac{1}{2}}\\
&\le C\int_{\{y\in \R^{2}:\;|y|\geq E\lambda_jG^{-1}_{a_i}(a_j)\}}\left(\frac{1}{1+|y|^2}\right)^{\frac{3}{2}}\left[\frac{\lambda_j}{\lambda_i}+\lambda_i\lambda_jG^{-2}_{a_i} (a_j )\right]^{-\frac{1}{2}}\\
&\le C\left[\frac{\lambda_j}{\lambda_i}+\lambda_i\lambda_jG^{-2}_{a_i} (a_j )\right]^{\frac{-1}{2}}\left(\frac{1}{\lambda_jG^{-1}_{a_i} (a_j )}\right).
\end{split}
\end{equation*}

\noindent
This implies
 \begin{equation}\label{i24fin}
 I^2_4=O\left(\varepsilon^2_{ij}\right).
 \end{equation}
Thus, combining \eqref{i4ine} and \eqref{i24fin}, we have that if \;$d_g(a_i,a_j)<2\d,$\; then 
\begin{equation}\label{i4pcclose}
I_4=O\left(\varepsilon^2_{ij}\right).
\end{equation}

\noindent
Now, using \eqref{part2i4} and \eqref{i4pcclose}, we infer that  in case \;$\varepsilon_{i, j}^{-2}\simeq \l_i\l_jG_{a_i}^{-2}(a_j),$\; 
\begin{equation}\label{i4pc}
I_4=O\left(\varepsilon_{i, j}^2\frac{1}{\d}\right).
\end{equation}
\noindent
Finally combining \eqref{part1i4}-\eqref{i4pc}, we get
\begin{equation}\label{I_44}
    I_4=O\left(\varepsilon^2_{ij}\frac{1}{\delta}\right).
\end{equation} 

\noindent
Using \eqref{J_11}-\eqref{I_33}, and \eqref{I_44}, we obtain the following for \;$J_1$\; (see \eqref{J1+J2})  
\begin{equation}\label{J_111}
J_1=c_0^4\left[1+O\left(\d+\frac{1}{\lambda^2_i\delta^2}\right)\right]\left[c_1\varepsilon_{ij}\left(1+o_{\varepsilon_{ij}}(1)+O(\varepsilon_{ij} \delta^{-1})\right)\right].
\end{equation}
\noindent
Thus, using \eqref{J1+J2}, \eqref{J_22}, and \eqref{J_111}, we arrive to 
\begin{equation*}
\begin{split}
\oint_{\partial M} u^3_{a_j,\lambda_j}u_{a_i,\lambda_i}\;dS_g=&c_0^4\left[1+O\left(\delta+\frac{1}{\lambda^2_i\delta^2}\right)\right]\left[c_1\varepsilon_{ij}\left(1+o_{\varepsilon_{ij}}(1)+O(\varepsilon_{ij} \delta^{-1})\right)\right]\\+&O\left(\varepsilon^2_{ij}\frac{1}{\delta^4}\right).
\end{split}
\end{equation*}

\noindent
Therefore, we obtain
\begin{equation}\label{eqfinal}
\begin{split}
\oint_{\partial M} u^3_{a_j,\lambda_j}u_{a_i,\lambda_i}\;dS_g=&c_0^4c_1\varepsilon_{ij}\left[\left(1+O\left(\delta+\frac{1}{\lambda^2_i\delta^2}\right)\right)\left(1+o_{\varepsilon_{ij}}(1)+O(\varepsilon_{ij} \delta^{-1})\right)\right]\\+&O\left(\varepsilon_{ij}^2\frac{1}{\delta^4}\right).
\end{split}
\end{equation}
Hence, recalling \;$\epsilon_{ji}=\oint_{\partial M} u^3_{a_j,\lambda_j}u_{a_i,\lambda_i}\;dS_g,$\;  then the result follows from \eqref{eqfinal}.
\end{pf}
\vspace{8pt}

\noindent
Clearly switching the index \;$i$\; and \;$j$\; in Lemma \ref{interact3}, we have the following corollary which is equivalent to Lemma \ref{interact3}. We decide to present the following corollary, because its form suits more our presentation of the Barycenter technique of Bahri-Coron\cite{bc} which follows the work \cite{martndia2} as done in \cite{nss}.
\begin{cor}\label{interact4}
Assuming that \;$\theta>0$\; is small and \;$\mu_0>0$\; is small then \;$\forall$\;\;$ a_i, a_j\in \partial M,$\;\;$\forall\; 0<2\delta<\delta_0,$\; and \;$\forall$\;\;$0<\frac{1}{\l_i}\le\frac{1}{\l_j}\leq \theta\delta$\; such that \;$\varepsilon_{ij}\leq \mu_0,$\; we have
\[\epsilon_{ij}=c_0^4c_1\varepsilon_{ij}\left[\left(1+O\left(\delta+\frac{1}{\lambda^2_j\delta^2}\right)\right)\left(1+o_{\varepsilon_{ij}}(1)+O(\varepsilon_{ij} \delta^{-1})\right)\right]+O\left(\varepsilon_{ij}^2\frac{1}{\delta^4}\right),\]
where \;$\delta_0$\; is as in \eqref{delta0}
\end{cor}
\vspace{8pt}

\noindent
We now show some sharp high-order inter-action estimates that are required for the application of the algebraic topological argument of Bahri-Coron\cite{bc} for existence. We begin with the balanced high-order inter-action estimate shown below.
\begin{lem}\label{interact5}
Assuming that \;$\theta>0$\; is small and \;$\mu_0>0$\; is small then \;$\forall$\;\;$ a_i, a_j\in \partial M$,\;\;$\forall 0<2\delta<\delta_0,$\; and \;$\forall$\;\;$0<\frac{1}{\l_i},$\;\;$\frac{1}{\l_j}\leq \theta\delta$\; such that \;$\varepsilon_{ij}\leq \mu_0,$\; we have
\[\oint_{\partial M} u^{2}_{a_i,\lambda_i} u^{2}_{a_j,\lambda_j}\;dS_g=O\left(\frac{\varepsilon^{2}_{ij}}{\delta^4}\log\left(\varepsilon^{-1}_{ij}\delta^{-1}\right)\right).\]
\end{lem}
\begin{pf}
By symmetry, we can assume without loss of generality (w.l.o.g) that\;$\lambda_j\le\lambda_i.$\; Thus we have \\
1) Either \;$\varepsilon^{-2}_{ij}\sim \lambda_i\lambda_jG^{-2}_{a_i}(a_j).$\;
\vspace{4pt}

\noindent
2) Or \;$\varepsilon^{-2}_{ij}\sim\frac{\lambda_i}{\lambda_j}.$\;
\vspace{6pt}

\noindent
Now, if \;$d_g(a_i,a_j)\geq2\delta,$\; then we have 
\begin{equation}\label{i1}
\begin{split}
I&:=\oint_{\partial M} u^{2}_{a_i,\lambda_i} u^{2}_{a_j,\lambda_j}\;dS_g\\
&\le C\oint_{\hat{B}(a_i,\delta)}\left(\frac{\lambda_i}{1+\lambda^2_id^{2}_g(a_i,x)}\right) \left(\frac{\lambda_j}{1+\lambda^2_jG^{-2}_{a_j}(x)}\right)\;dS_g\\
&+\frac{C}{\lambda_i\delta^2}\oint_{\hat{B}(a_j,\delta)}\left(\frac{\lambda_j}{(1+\lambda^2_jd^{2}_g(a_j,x)}\right)\;dS_g+\frac{C}{\lambda_i\lambda_j\delta^4}\\
&\le\underbrace{C\int_{\hat{B}_{\delta\lambda_i}(0)}\frac{1}{1+|y|^2}\left(\frac{1}{\frac{\lambda_i}{\lambda_j}+\lambda_i\lambda_jG^{-2}_{a_j}\left(\psi_{a_i}\left(\frac{y}{\lambda_i}\right)\right)}\right)}_{\mbox{$I_1$}}\\
&+\frac{C}{\lambda_i\lambda_j\delta^2}\int_{\hat{B}_{\delta\lambda_j}(0)}\left(\frac{1}{1+|y|^2}\right)+\frac{C}{\lambda_i\lambda_j\delta^4}\\&
\le C I_1+\frac{C}{ \lambda_i\lambda_j \delta^2}\left[\log\left(\lambda_j\delta\right)+C\right]+\frac{C}{\lambda_i\lambda_j\delta^4}\\&
\le C I_1+\frac{C}{\lambda_i\lambda_j\delta^4}\log\left(\lambda_j\right).
\end{split}
\end{equation}
\noindent
Now, we estimate \;$I_1$\; as follows\\\\
If \;$\varepsilon^{-2}_{ij}\sim\frac{\lambda_i}{\lambda_j},$\; then we get
\[I_1\le C \varepsilon^2_{ij}\left[\log\left(\lambda_i\delta\right)+C\right].\]
So, for \;$I$\; we have
\begin{equation*}
\begin{split}
I&\le C \varepsilon^2_{ij}\left[\log\left(\lambda_i\delta\right)+C\right]+\frac{C}{\lambda_i\lambda_j\delta^4}\log\left(\lambda_i\lambda_j\right)
\\&\le\frac{C}{\delta^4}\varepsilon^{2}_{ij}\log\left(\varepsilon^{-2}_{ij}G^{2}_{a_i}(a_j)\right)\\
&=O\left(\frac{\varepsilon^2_{ij}\log\left(\varepsilon^{-1}_{ij}\delta^{-1}\right)}{\delta^4}\right).
\end{split}
\end{equation*}

\noindent
If \;$\varepsilon^{-2}_{ij}\sim \lambda_i\lambda_jG^{-2}_{a_i}(a_j),$\; then we get
\begin{equation*}
I_1\le \frac{C}{\lambda_i\lambda_j\delta^2}\left[\log\left(\lambda_i\delta\right)+C\right].
\end{equation*}
\noindent
So, for \;$I$\; we have
\begin{equation*}
I\le\frac{C}{\lambda_i\lambda_j\delta^4}\left[\log\left(\lambda_i\lambda_j\right)\right].
\end{equation*}
This implies 
\begin{equation*}
I\le \frac{C}{\delta^4}\varepsilon^{2}_{ij}\log\left(\varepsilon^{-2}_{ij}G^{2}_{a_i}(a_j)\right).
\end{equation*} 

\noindent
Hence, for \;$d_g(a_i,a_j)\geq 2\d,$\; we obtain
\begin{equation}\label{estfar}
I=O\left(\frac{\varepsilon^{2}_{ij}\log\left(\varepsilon^{-1}_{ij}\delta^{-1}\right)}{\delta^4}\right).
\end{equation}

\noindent
On the other hand, arguing as above, if \;$d_g(a_i,a_j)<2\delta,$\; then we have also 
\begin{equation*}
\begin{split}
I&\le I_1+\frac{C}{\lambda_i\lambda_j\delta^2}\left[\log\left(\lambda_j\right)\right]+\frac{C}{\lambda_i\lambda_j\delta^4}\\&
\le I_1+\frac{C}{\lambda_i\lambda_j\delta^4}\log\left(\lambda_i\lambda_j\right),
\end{split}
\end{equation*}
where \;$I_1$\; is as in \eqref{i1}. Thus, if \;$\varepsilon_{i,j}^{-2}\simeq \frac{\l_i}{\l_j},$\; then
\begin{equation*}
  \begin{split}
      I\le I_1+\frac{C}{\delta^4}\left(\frac{\lambda_j}{\lambda_i}\right)\frac{1}{\lambda_j}\left[\log(\frac{\lambda_i}{\lambda_j})+\log(\lambda^{2}_j)\right].
  \end{split}  
\end{equation*}
This implies
\begin{equation*}
    I\leq  I_1+\frac{C}{\d^4}\varepsilon_{ij}^2\log (\varepsilon_{ij}^{-1}).
\end{equation*}
\vspace{5pt}
\noindent 
Next, if \;$\varepsilon_{i,j}^{-2}\simeq \l_i\l_j G^{-2}_{a_i}(a_j),$\; then we get
\begin{equation*}
\begin{split}
 I&\le I_1+\frac{1}{\lambda_i\lambda_j\delta^4G^{-2}_{a_i}(a_j)}\left[\log\left(\lambda_i\lambda_jG^{-2}_{a_i}(a_j)\right)+\log\left(G^{2}_{a_i}(a_j)\right)\right]G^{-2}_{a_i}(a_j)\\
 &\le I_1+\frac{C}{\delta^4}\varepsilon_{ij}^2\log\left(\varepsilon_{ij}^{-1}\right).
\end{split}
\end{equation*}

\noindent
Now, in order to proceed, we will estimate \;$I_1.$\; To do so, we begin by defining the following sets:
\begin{equation*}
\begin{split}
&A_1=\left\{y\in \R^{2}:\;|y|\le \epsilon\lambda_i \sqrt{G^{-2}_{a_j}(a_i)+\frac{1}{\lambda^2_j}}\right\},\\&
A_2=\left\{y\in \R^{2}:\;\epsilon\lambda_i \sqrt{G^{-2}_{a_j}(a_i)+\frac{1}{\lambda^2_j}}\le |y|\le E \lambda_i\sqrt{G^{-2}_{a_j}(a_i)+\frac{1}{\lambda^2_j}}\right\},\\&
A_3=\left\{y\in \R^{2}:\;E\lambda_i \sqrt{G^{-2}_{a_j}(a_i)+\frac{1}{\lambda^2_j}}\le |y|\le 4\lambda_i\delta\right\},
\end{split}
\end{equation*}
with \;$0<\epsilon<E<\infty$. Clearly by  the definition of \;$I_1$ (see \eqref{i1}), we have
\[I_1\le \int_{A_1}L_{ij}+\int_{A_2}L_{ij}+\int_{A_3}L_{ij},\]
where 
\[L_{ij}=\left(\frac{1}{1+|y|^2}\right)\left(\frac{1}{\frac{\lambda_i}{\lambda_j}+\lambda_i\lambda_jG^{-2}_{a_j}\left(\psi_{a_i}\left(\frac{y}{\lambda_i})\right)\right)}\right).\]

\noindent
For \;$\int_{A_1}L_{ij},$\; we have 
\begin{equation*}\begin{split}
\int_{A_1}L_{ij} &\le C \varepsilon^{2}_{ij}\int_{A_1}\left(\frac{1}{1+|y|^2}\right)\\
&\le C\varepsilon^{2}_{ij} \log\left(\sqrt{\frac{\lambda_i}{\lambda_j}}\sqrt{\lambda_i\lambda_jG^{-2}_{a_j}(a_i)+\frac{\lambda_i}{\lambda_j}}\right)\\&
\le C\varepsilon^{2}_{ij} \log\left(\varepsilon^{-1}_{ij}\right).
\end{split}
\end{equation*}

\noindent
For \;$\int_{A_2}L_{i,j},$\; we have
\begin{equation*}
\begin{split}
\int_{A_2}L_{ij}&\le C \left(\frac{1}{\left(\frac{\lambda_i}{\lambda_j}\right)^2+\lambda^2_iG^{-2}_{a_j}(a_i)}\right) \int_{A_2}\left(\frac{1}{\frac{\lambda_i}{\lambda_j}+\lambda_i\lambda_jG^{-2}_{a_j}\left(\psi_{a_i}\left(\frac{y}{\lambda_i}\right)\right)}\right)\\&
\le C \left(\frac{\lambda_j}{\lambda_i}\right) \varepsilon^{2}_{ij}\int_{\{y\in\R^{2}:\;|y|\leq E\l_i\sqrt{G_{a_j}^{-2}(a_i)+\frac{1}{\l_j^2}}\}}\left(\frac{1}{\frac{\lambda_i}{\lambda_j}+\frac{\lambda_j}{\lambda_i}\left|\lambda_i\psi^{-1}_{a_j}\circ\psi_{a_i}\left(\frac{y}{\lambda_i}\right)\right|^2}\right)\\
&\le C \left(\frac{\lambda_j}{\lambda_i}\right)\varepsilon^{2}_{ij}\int_{\{y\in\R^{2}:\;|y|\leq \bar E\l_i\sqrt{G_{a_j}^{-2}(a_i)+\frac{1}{\l_j^2}}\}}\left(\frac{1}{\frac{\lambda_i}{\lambda_j}+\frac{\lambda_j}{\lambda_i}\left|y\right|^2}\right)\\
&\le C \left(\frac{\lambda_j}{\lambda_i}\right)^{2}\left(\frac{\lambda_j}{\lambda_i}\right)^{-2}\varepsilon^{2}_{ij}\int_{\{y\in\R^{2}:\;|y|\leq \bar E\l_j\sqrt{G_{a_j}^{-2}(a_i)+\frac{1}{\l_j^2}}\}}\left(\frac{1}{1+|y|^2}\right)\\&
\le C\varepsilon^{2}_{ij}\log\left(\varepsilon^{-1}_{ij}\right).
\end{split}
\end{equation*}

\noindent
For \;$\int_{A_3}L_{i,j},$\; we have 
\begin{equation*}
\begin{split}
\int_{A_3}L_{ij}&\le\int_{A_3}\left(\frac{1}{1+|y|^2}\right)\left(\frac{1}{\frac{\lambda_i}{\lambda_j}+\frac{\lambda_j}{\lambda_i}|y|^2}\right)\\&
\le C \left(\frac{\lambda_i}{\lambda_j}\right)\int_{A_3}\frac{1}{|y|^4}\\&
\le C \left(\frac{\lambda_i}{\lambda_j}\right)\frac{1}{\left(\lambda^{2}_iG^{-2}_{a_j}(a_i)+\left(\frac{\lambda_i}{\lambda_j}\right)^2\right)}\\&
\le C \varepsilon^{2}_{ij}.
\end{split}
\end{equation*}
\noindent
Therefore, we have  
\begin{equation*}
    \begin{split}
      I_1\leq C\varepsilon^{2}_{ij}\log{\varepsilon^{-1}_{ij}}.  
    \end{split}
\end{equation*}

\noindent
This implies for \;$d_g(a_i,a_j)<2\d,$\;  we have 
\begin{equation*}
    \begin{split}
       I=O\left(\frac{\varepsilon^{2}_{ij}}{\delta^4}\log\left(\varepsilon^{-1}_{ij}\right)\right). 
    \end{split}
\end{equation*}
 
\noindent
Hence, combining with the estimate for \;$d_g(a_i,a_j)\geq 2\d$\; (see \eqref{estfar}), we have 
\[\oint_{\partial M} u^2_{a_i,\lambda_i} u^2_{a_j,\lambda_j}\;dS_g=O\left(\frac{\varepsilon^{2}_{ij}}{\delta^4}\log\left(\varepsilon^{-1}_{ij}\delta^{-1}\right)\right).\]
\end{pf}\\

\noindent
\vspace{4pt}
Finally, we establish a sharp unbalanced high-order inter-action estimate that is required for the application of the Barycenter technique of Bahri-Coron\cite{bc} for existence. 

\noindent
\begin{lem}\label{interact6}
Assuming that \;$\theta>0$\; is small and \;$\mu_0>0$\; is small, then \;$\forall$\;\;$ a_i, a_j\in \partial M$,\;\;$\forall 0<2\delta<\delta_0,$\; and \;$\forall$\;\;$0<\frac{1}{\l_i}\le\frac{1}{\l_j}\leq \theta\delta$\; such that \;$\varepsilon_{ij}\leq \mu_0,$\; we have
\[\oint_{\partial M} u^{\alpha}_{a_i,\lambda_i} u^{\beta}_{a_j,\lambda_j}\;dS_g=O\left(\frac{\varepsilon_{ij}^{\beta}}{\d^4}\right).\]
where \;$\delta_0$\; is as in \eqref{delta0}, \;\(\alpha+\beta=4,\)\; and \;\(\alpha>2>\beta>1.\)
\end{lem}
\begin{pf}
Let \;$\hat{\alpha}=\frac{1}{2}\alpha$\; and \;$\hat{\beta}=\frac{1}{2}\beta.$\; Then we have \;$\hat{\alpha}+\hat{\beta}=2.$\; Now, since \;$\lambda_j\le \lambda_i,$\; then for \;$\varepsilon_{ij}\sim 0$\; we have\\
1) Either \;$\varepsilon^{-2}_{ij}\sim \lambda_i\lambda_jG^{-2}_{a_i}(a_j).$\;
\vspace{4pt}

\noindent
2) Or \;$\varepsilon^{-2}_{ij}\sim\frac{\lambda_i}{\lambda_j}.$\;
\vspace{4pt}

\noindent
To continue, we write
\[\oint_{\partial M} u^{\alpha}_{a_i,\lambda_i} u^{\beta}_{a_j,\lambda_j}\;dS_g=\underbrace{\oint_{\hat{B}(a_i,\delta)} u^{\alpha}_{a_i, \lambda_i} u^{\beta}_{a_j,\lambda_j}\;dS_g}_{\mbox{$I_1$}}+\underbrace{\oint_{{\partial M}-\hat{B}(a_i,\delta)} u^{\alpha}_{a_i, \lambda_i} u^{\beta}_{a_j, \lambda_j}\;dS_g}_{\mbox{$I_2$}}\]
\noindent
and estimate \;$I_1$\; and \;$I_2.$\; For \;$I_2,$\; we have 
\begin{equation*}\begin{split}
I_2&= \oint_{\left({\partial M}-\hat{B}(a_i,\delta)\right)\cap B(a_j,\delta)} u^{\alpha}_{a_i,\lambda_i} u^{\beta}_{a_j,\lambda_j}\;dS_g+\oint_{{\partial M}-\left(\hat{B}(a_i,\delta)\cup \hat{B}(a_j,\delta)\right)} u^{\alpha}_{a_i,\lambda_i} u^{\beta}_{a_j,\lambda_j}\;dS_g\\
&\le C\oint_{\left({\partial M}-\hat{B}(a_i,\delta)\right)\cap \hat{B}(a_j,\delta)} \left(\frac{\lambda_i}{1+\lambda^2_iG^{-2}_{a_i}(x)}\right)^{\hat{\alpha}}\left(\frac{\lambda_j}{1+\lambda^2_jd^{2}_g(a_j,x)}\right)^{\hat{\beta}}\;dS_g\\
&+C\oint_{{\partial M}-\left(\hat{B}(a_i,\delta)\cup \hat{B}(a_j,\delta)\right)}\left(\frac{\lambda_i}{1+\lambda^2_iG^{-2}_{a_i}(x)}\right)^{\hat{\alpha}}\left(\frac{\lambda_j}{1+\lambda^2_jG^{-2}_{a_j}(x)}\right)^{\hat{\beta}}\;dS_g\\
&\le \frac{C}{\lambda^{\hat{\alpha}}_i\lambda^{2-\hat{\beta}}_j\delta^{\alpha}}\int_{\hat{B}_{\delta\lambda_j}(0)}\left(\frac{1}{1+|y|^2}\right)^{\hat{\beta}}+\frac{C}{\lambda^{\hat{\alpha}}_i\lambda^{\hat{\beta}}_j\delta^{4}}\\&
\le \frac{C}{\lambda^{\hat{\alpha}}_i\lambda^{2-\hat{\beta}}_j\delta^{\alpha}}\left(\frac{1}{\lambda_j \delta}\right)^{2\hat{\beta}-2}+\frac{C}{\lambda^{\hat{\alpha}}_i\lambda^{\hat{\beta}}_j\delta^{4}}.
\end{split}
\end{equation*}

\noindent
Thus, we have for \;$I_2$\;
\begin{equation}\label{esti2int}
I_2\le \frac{C}{\lambda^{\hat{\alpha}}_i\lambda^{\hat{\beta}}_j\delta^{4}}.
\end{equation}
\noindent
Next, for \;$I_1$\; we have 
\begin{equation*}
\begin{split}
I_1&=\oint_{\hat{B}(a_i,\delta)}\left(\frac{\lambda_i}{(1+\lambda^2_id^{2}_g(a_i,x)}\right)^{\hat{\alpha}}\left(\frac{\lambda_j}{1+\lambda^2_jG^{-2}_{a_j}(x)}\right)^{\hat{\beta}}\;dS_g\\
&= \int_{\hat{B}_{\delta\lambda_i}(0)} \left(\frac{1}{1+|y|^2}\right)^{\hat{\alpha}}\left[\frac{1}{\frac{\lambda_i}{\lambda_j}+\lambda_i\lambda_jG^{-2}_{a_j}\left(\psi_{a_i}\left(\frac{y}{\lambda_i}\right)\right)}\right]^{\hat{\beta}}.
\end{split}\end{equation*}
\noindent
Thus, if \;$\varepsilon^{-2}_{ij}\sim \frac{\lambda_i}{\lambda_j},$\; then 
\begin{equation*}
\begin{split}
I_1&\le C \varepsilon^{2\hat{\beta}}_{ij}\left[\left(\frac{1}{\lambda_i\delta}\right)^{2\hat{\alpha}-2}+C\right]\\&\le C \varepsilon^{\beta}_{ij}.
\end{split}
\end{equation*}
\noindent
If \;$\varepsilon^{-2}_{ij}\sim \lambda_i \lambda_j G^{-2}_{a_i}(a_j)$\; and \;$d_g(a_i,a_j)\geq2\delta,$\; then we have 
\begin{equation*}
\begin{split}
I_1&\le C \left(\frac{1}{\lambda_i\lambda_j\delta^2}\right)^{\hat{\beta}}\left[\left(\frac{1}{\lambda_i\delta}\right)^{2\hat{\alpha}-2}+C\right]\\&
\le C \frac{1}{\delta^2}\left(\frac{1}{\lambda_i\lambda_j}\right)^{\hat{\beta}}\le C \frac{1}{\delta^2}\left[\left(\frac{1}{\lambda_i\lambda_jG^{-2}_{a_i}(a_j)}\right)^{\frac{1}{2}}\right]^{\beta}
\\&\le C \frac{1}{\delta^2}\varepsilon^{\beta}_{ij}.
\end{split}
\end{equation*}
\noindent
Now, if \;$\varepsilon^{-2}_{ij}\sim \lambda_i \lambda_j G^{-2}_{a_i}(a_j)$\; and \;$d_g(a_i,a_j)<2\delta,$\; then we get
\[I_1\le C\int_{\hat{B}_{\delta\lambda_i}(0)} \left(\frac{1}{1+|y|^2}\right)^{\hat{\alpha}}\left[\frac{1}{\frac{\lambda_i}{\lambda_j}+\lambda_i\lambda_j\left|\psi^{-1}_{a_j}\circ\psi_{a_i}\left(\frac{y}{\lambda_i}\right)\right|^{2}}\right]^{\hat{\beta}}.\]
Next, we define 
\[B=\left\{y\in\R^{2}:\;\frac{1}{2}d_g(a_i,a_j)\le \frac{|y|}{\lambda_i}\le2d_g(a_i,a_j)\right\}\]
and have
\begin{equation*}
\begin{split}
I_1&\le C\int_{B}\left(\frac{1}{1+|y|^2}\right)^{\hat{\alpha}}\left[\frac{1}{\frac{\lambda_i}{\lambda_j}+\lambda_i\lambda_j\left|\psi^{-1}_{a_j}\circ\psi_{a_i}\left(\frac{y}{\lambda_i}\right)\right|^{2}}\right]^{\hat{\beta}}\\
&+C\int_{\hat{B}_{\delta\lambda_i}(0)-B} \left(\frac{1}{1+|y|^2}\right)^{\hat{\alpha}}\left[\frac{1}{\frac{\lambda_i}{\lambda_j}+\lambda_i\lambda_j\left|\psi^{-1}_{a_j}\circ\psi_{a_i}\left(\frac{y}{\lambda_i}\right)\right|^{2}}\right]^{\hat{\beta}}.
\end{split}
\end{equation*}
\noindent
For the second term, we have
\begin{equation*}
\begin{split}
\int_{\hat{B}_{\delta\lambda_i}(0)-B} \left(\frac{1}{1+|y|^2}\right)^{\hat{\alpha}}\left[\frac{1}{\frac{\lambda_i}{\lambda_j}+\lambda_i\lambda_j\left|\psi^{-1}_{a_j}\circ\psi_{a_i}\left(\frac{y}{\lambda_i}\right)\right|^{2}}\right]^{\hat{\beta}}
&\le C \varepsilon^{\beta}_{ij}\left[\left(\frac{1}{\lambda_i\delta}\right)^{\alpha-2}+C\right]\\&\le C \varepsilon^{\beta}_{ij}.
\end{split}
\end{equation*}
\noindent
For the first term, we have
\begin{equation*}
\begin{split}
&\int_{B}\left(\frac{1}{1+|y|^2}\right)^{\hat{\alpha}}\left[\frac{1}{\frac{\lambda_i}{\lambda_j}+\lambda_i\lambda_j\left|\psi^{-1}_{a_j}\circ\psi_{a_i}\left(\frac{y}{\lambda_i}\right)\right|^{2}}\right]^{\hat{\beta}}\\
&\le C\left(\frac{1}{1+\lambda^2_id^{2}_g(a_i,a_j)}\right)^{\hat{\alpha}}\int_{\{y\in\R^{2}:\;|y|\le 2\lambda_id_g(a_i,a_j)\}}\left[\frac{1}{\frac{\lambda_i}{\lambda_j}+\frac{\lambda_j}{\lambda_i}\left|\lambda_i\psi^{-1}_{a_j}\circ\psi_{a_i}\left(\frac{y}{\lambda_i}\right)\right|^2}\right]^{\hat{\beta}}\\
&\le C\left(\frac{1}{1+\lambda^2_id^{2}_g(a_i,a_j)}\right)^{\hat{\alpha}}\int_{\{z\in\R^{2}:\;|z|\le 4\lambda_id_g(a_i,a_j)\}}\left[\frac{1}{\frac{\lambda_i}{\lambda_j}+\frac{\lambda_j}{\lambda_i}|z|^2}\right]^{\hat{\beta}}\\
&\le C\left(\frac{1}{\frac{\lambda_j}{\lambda_i}+\lambda_i\lambda_jd_g(a_i,a_j)^2}\right)^{\frac{\alpha}{2}}\int_{\{z\in\R^{2}:\;|z|\le 4\lambda_jd_g(a_i,a_j)\}}\left[\frac{1}{1+|z|^2}\right]^{\hat{\beta}}.
\end{split}
\end{equation*}
\noindent
If \;$\lambda_jd_g(a_i,a_j)$\; is bounded, then we get
\begin{equation*}
\begin{split}
I_1&\le C\left(\frac{1}{\frac{\lambda_j}{\lambda_i}+\lambda_i\lambda_jd^{2}_g(a_i,a_j)}\right)^{\frac{\alpha}{2}}\\&\le C\varepsilon^{\beta}_{ij}. 
\end{split}
\end{equation*}
\noindent
If \;$\lambda_jd_g(a_i,a_j)$\; is unbounded, then we get
\begin{equation*}
\begin{split}
I_1&\le C\left(\frac{1}{\frac{\lambda_j}{\lambda_i}+\lambda_i\lambda_jd^{2}_g(a_i,a_j)}\right)^{\frac{\alpha}{2}}\left(\lambda_jd^{2}_g(a_i,a_j)\right)^{2-2\hat{\beta}}\\
&\le C \left(\frac{1}{1+\lambda^2_id^{2}_g(a_i,a_j)}\right)^{\hat{\alpha}+\hat{\beta}-1}\left(\frac{\lambda_i}{\lambda_j}\right)^{\hat{\beta}}\\
&\le C \left(\frac{1}{\frac{\lambda_j}{\lambda_i}+\lambda_i\lambda_jd^{2}_g(a_i,a_j)}\right)^{\hat{\beta}}\left(\frac{1}{1+\lambda^2_id^{2}_g(a_i,a_j)}\right)^{\hat{\alpha}-1}\\&\le C \varepsilon^{\beta}_{ij}.
\end{split}
\end{equation*}
\noindent
Thus, we have for \;$I_1$\;
\begin{equation}\label{esti1f}
I_1\le \frac{C}{\delta^2}\varepsilon^{\beta}_{ij}.
\end{equation}
\noindent
On the other hand, using the estimate for \;$I_2$\; (see \eqref{esti2int}), we have
\begin{equation}\label{esti2f}
I_2=O\left(\frac{\varepsilon^{\beta}_{ij}}{\delta^4}\right).
\end{equation}
Hence, combining \eqref{esti1f} and \eqref{esti2f}, we have 
\[\oint_{\partial M} u^{\alpha}_{a_i,\lambda_i} u^{\beta}_{a_j,\lambda_j}\;dS_g=O\left(\frac{\varepsilon^{\beta}_{ij}}{\delta^4}\right).\]
\end{pf}
\vspace{4pt}
%
%
%
%
%
\section{Algebraic topological argument}\label{ATA}
In this Section, we present the algebraic topological argument for existence. We start by fixing some notations from algebraic topology.  
For  a topological space \;$Z$\; and \;$Y$\; a subspace of \;$Z,$\;\; $H_*(Z, Y)$\; stands for the relative homology with \;$\Z_2$\; coefficients of the topological pair \;$(Z, Y)$. For\;$f: (Z, Y)\longrightarrow (W, X)$ \; a map with \;$(Z, Y)$\; and \;$(W, X)$\; topological pairs, \;$ f_*$ denotes the induced map in relative homology.
\vspace{4pt}

\noindent
Furthermore, we discuss some algebraic topological tools needed for our application of the Barycenter technique of Bahri-Coron\cite{bc} for existence. We start with recalling the space of formal the  barycenter of \;$\partial M.$\; For \;$p\in \N^*$,\; the set of formal barycenters of \;$\partial M$\; of order \;$p$\; is defined as 
 \begin{equation}\label{eq:barytop}
B_{p}(\partial M)=\{\sum_{i=1}^{p}\alpha_i\d_{a_i}\;:\;a_i\in \partial M, \;\alpha_i\geq 0,\;\; i=1,\cdots, p,\;\,\sum_{i=1}^{p}\alpha_i=1\},\;\;\text{and}\;\;B_0(\partial M)=\emptyset,
\end{equation}
where \;$\delta_{a}$\; for \;$a\in \partial M$\; is the Dirac measure at \;$a.$\; Since \;$dim(\partial M)=2,$\; then we have the existence of \;$\Z_2$\; orientation classes  (see \cite{bc} and \cite{kk})
\begin{equation}\label{orientation_classes}
w_p\in H_{2p-1}(B_{p}(\partial M), B_{p-1}(\partial M)), \;\;\;\;\;p\in \N^*.
\end{equation}\\
Now to continue, we fix \;$\delta$\; small such that \;$0<2\d< \delta_0$ \;where \;$\delta_0$ \; is as in \eqref{delta0}. Moreover, we choose \;$\theta_0>0$\; and smalll. After this, we let \;$\l$\; varies such that \;$0<\frac{1}{\l}\leq \theta_0\delta$\; and associate for every \;$p\in \N^*$\; the map
$$
f_p(\l): B_p(\partial M)\longrightarrow H^{1}_{+}(M) 
$$
defined by the formula
$$
f_p(\l)(\sigma)=\sum_{i=1}^p\alpha_i u_{a_i,\l}, \;\;\;\;\sigma=\sum_{i
=1}^p\alpha_i\d_{a_i}\in B_p(\partial M),
$$
where \;$u_{a_i, \l}$\; is as in \eqref{ual} with \;$a$\; replaced by \;$a_i$\;.\\
\vspace{6pt}

\noindent
As in Proposition 3.1 in  \cite{martndia2} and Proposition 6.3 in \cite{nss}, using Corollary \ref{sharpenergy}, Corollary \ref{interact2}, Corollary \ref{interact4}, Lemma \ref{interact5}, and  Lemma \ref{interact6}, we have the following multiple-bubble estimate.
\begin{pro}\label{eq:baryest}
There exist \;$\bar C_0>0$\; and \;$\bar c_0>0$\; such that for every \;$p\in \N^*$, $p\geq 2$ and every \;$0<\varepsilon\leq \varepsilon_0$, there exists \;$\l_p:=\l_p(\varepsilon)$ such that for every \;$\l\geq \l_p$ and for every $\sigma=\sum_{i=1}^p\alpha_i\delta_{a_i}\in B_p(\partial M)$, we have
\begin{enumerate}
 \item
If \;$\sum_{i\neq j}\varepsilon_{i, j}> \varepsilon$\; or there exist \;$i_0\neq j_0$\; such that \;$\frac{\alpha_{i_0}}{\alpha_{j_0}}>\nu_0$, then
 $$
J_q(f_p(\l)(\sigma))\leq p^{\frac{1}{2}}\mathcal{S}.
$$ 
 \item
If \;$\sum_{i\neq j}\varepsilon_{i, j}\leq \varepsilon$\; and for every \;$i\neq j$\; we have \;$\frac{\alpha_{i}}{\alpha_j}\leq\nu_0$, then
$$
J_q(f_p(\l)(\sigma))\leq p^{\frac{1}{2}}\mathcal{S}\left(1+ \frac{\bar C_0}{\l}-\bar c_{0}\frac{(p-1)}{\l}\right),
$$
where \;$\mathcal{S}$\; is as in \eqref{S}, \;$\varepsilon_{ij}$\; is as in \eqref{varepij}, \;$\l_i=\l_j=\l,$\;\;$\varepsilon_0$\; is as in \eqref{varepsilon0}, and \;$\nu_0$\; is as in \eqref{nu0}.
\end{enumerate}

\end{pro}
\vspace{6pt}

\noindent
As in Lemma 4.2  in \cite{martndia2} and Lemma 6.4 in \cite{nss}, we have the selection map \;$s_1$ (see \eqref{eq:mini}), Lemma \ref{deform} and Corollary \ref{sharpenergy} imply the following topological result.
\begin{lem}\label{eq:nontrivialf1}
Assuming that \;$J_q$\; has no critical points, then there exists \;$\bar \l_1>0$\;  such that for every \;$\l\geq \bar \l_1$, we have
$$
f_1(\l)\; : \;(B_1(\partial M),\; B_0(\partial M))\longrightarrow (W_1, \;W_0)
$$
is well defined and satisfies
$$
(f_1(\l))_*(w_1)\neq 0\;\;\text{in}\;\;H_{2}(W_1, \;W_0).
$$
\end{lem}
\vspace{6pt}

\noindent
As in Lemma 4.3  in \cite{martndia2} and Lemma 6.5 in \cite{nss}, we have the selection map \;$s_p$ (see \eqref{eq:mini}), Lemma \ref{deform} and Proposition \ref{eq:baryest} imply the following recursive topological result.
\begin{lem}\label{eq:nontrivialrecursive}
Assuming that \;$J_q$\; has no critical points, then  there exists \;$\bar \l_p>0$\; such that for every \;$\l\geq \bar\l_p$, we have 
$$
f_{p+1}(\l): (B_{p+1}(\partial M),\; B_{p}(\partial M))\longrightarrow (W_{p+1}, \;W_{p})
$$
and 
$$
f_p(\l): (B_p(\partial M), \;B_{p-1}(\partial M))\longrightarrow (W_p, \; W_{p-1})
$$
are well defined and satisfy
$$(f_p(\l))_*(w_p)\neq 0\;\; \text{in}\;\; \;\;H_{2p-1}(W_p, \;W_{p-1})$$ implies
$$(f_{p+1}(\l))_*(w_{p+1})\neq 0\;\; \text{in} \;\;H_{2(p+1)-1}(W_{p+1}, \;W_{p}).$$
\end{lem}
\vspace{6pt}

\noindent
Finally, as in Corollary 3.3 in \cite{martndia2} and Lemma 6.6 in \cite{nss}, we clearly have that Proposition \ref{eq:baryest}  implies the following result. 
\begin{lem}\label{eq:largep}
Setting \;$$\bar p_{0}:=[1+\frac{\bar C_{0}}{\bar c_{0}} ]+2$$ with \;$\bar C_0$\; and \;$\bar c_0$\; as in Proposition \ref{eq:baryest} and recalling \eqref{dfwp}, we have  there exists \;$\hat \l_{\bar p_0}>0$\; such that \;$\forall\l\geq\hat \l_{\bar p_0}$, 
$$
f_{\bar p_{0}}(\l)(B_{\bar p_{0}}(\partial M))\subset W_{{\bar p}_0-1}.
$$
\end{lem}

\noindent
\begin{pfn} {of Theorem \ref{thm1}} \\
As in \cite{martndia2} and \cite{nss}, the theorem follows by a contradiction argument from Lemma \ref{eq:nontrivialf1} - Lemma \ref{eq:largep}.
\end{pfn}

\section{Appendix}\label{APP}
In this Section, using the explicit expression of \;$\delta_{0,\lambda}$\; (see \eqref{bubble-M3}) or Lemma A-1 in \cite{nss}, we have the following technical estimates.
\begin{lem}\label{Appendix}
Recalling the definition of \;$\delta_{0,\lambda}$\; see \eqref{bubble-M3}, and setting \;$x=(\bar{x},x_3)$\; with \;$\bar{x}=(x_1,x_2),$\; we have on \;$\bar{\R}^{3}_{+}$\;
\end{lem}
 \begin{equation*}
     \begin{split}
         &\;\;\delta_{0,\lambda}(x)=O\left(\left(\frac{\lambda}{1+\lambda^{2}|x|^{2}}\right)^{\frac{1}{2}}\right),\\
         &\;\;\partial_{x_3}\delta_{0,\lambda}(x)=O\left(\left(\frac{\lambda}{1+\lambda^{2}|x|^{2}}\right)^{\frac{1}{2}}\left(\frac{\lambda^2}{1+\lambda^{2}|x|^{2}}\right)^{\frac{1}{2}}\right),\\
         &\;\;\nabla_{\bar{x}}\delta_{0,\lambda}(x)=O\left(\left(\frac{\lambda}{1+\lambda^{2}|x|^{2}}\right)^{\frac{1}{2}}\left(\frac{\lambda^2}{1+\lambda^{2}|x|^{2}}\right)^{\frac{1}{2}}\right),\\
         &\;\;\nabla^{2}_{\bar{x}}\delta_{0,\lambda}(x)=O\left(\left(\frac{\lambda}{1+\lambda^{2}|x|^{2}}\right)^{\frac{1}{2}}\left(\frac{\lambda^2}{1+\lambda^{2}|x|^{2}}\right)\right).
     \end{split}
 \end{equation*}
 
\end{document}